\documentclass[10pt]{article}
\usepackage{amsfonts, amsmath, amssymb, latexsym}

\numberwithin{equation}{section}
\newcommand{\ad}{\hbox{\rm ad\,}}
\newcommand{\Ann}{\hbox{\rm Ann\,}}

\newcommand{\Hom}{\hbox{\rm Hom\,}}
\newcommand{\Der}{\hbox {\rm Der}}
\newcommand{\eqdef}{\buildrel{\hbox{\small def}} \over =}

\newcommand{\qed}{\quad $~~\square$}
\newcommand{\keywords}[1]{\textbf{\textit{Keywords}} #1}
\newcommand{\subjectclassification}[1]{\textbf{\textit{ Subject Classification}} #1}
\newtheorem{Thm}[equation] {Theorem}
\newtheorem{Lem}[equation] {Lemma}
\newtheorem{Cor}[equation] {Corollary}

\newtheorem{Pro}[equation] {Proposition}

\newtheorem{Ex}[equation] {Example}

\begin{document}

\title{On Graded Lie Algebras of Characteristic Three \\
With Classical Reductive Null Component}
\author{T. B. Gregory \\ M. I. Kuznetsov\footnote{The second author gratefully acknowledges partial
support from the Russian Foundation of Basic Research Grants
\#{}02-01-00725 and \#{}18-01-00900.  He would also like to express his appreciation
for the hospitality of The Ohio State University, both at Columbus
and at Mansfield, and for the support of The Ohio State University
at Mansfield. \newline}}
\maketitle


\subjectclassification{180: Lie algebras and Lie superalgebras}

\setcounter{section}{-1}
\begin{section}{Introduction} \end{section}

\begin{abstract}We consider finite-dimensional irreducible transitive graded Lie
algebras $L = \sum_{i=-q}^rL_i$ over algebraically closed fields
 of characteristic three.  We assume that the null component $L_0$
 is classical and reductive.  The adjoint representation of $L$ on itself 
induces a representation of the commutator subalgebra $L_0'$ 
of the null component on the minus-one component $L_{-1}.$  
We show that if the depth $q$ of $L$ is greater than one,
 then this representation must be restricted.\end{abstract}

     Over algebraically closed fields ${\mathbf F}$ of characteristic $p>0,$ the classification of the finite-dimensional simple Lie algebras relies on the classification of the
finite-dimensional irreducible transitive graded Lie algebras
$L=\bigoplus_{i=-q}^r L_i$ of depth $q\geqq 1$ with classical
reductive null component $L_0.$
 We recall some of the progress that has
been made in the classification of such Lie algebras $L.$ In the
case in which $L_{-1}$ is not only irreducible but also restricted
as an $L_0$-module, such Lie algebras are described by the
Recognition Theorem of Kac \cite{K} for $p>5.$ (See  also
\cite{BGP}.)  In \cite{BG} it is shown that for $p>5,$ $L_{-1}$ is
necessarily a restricted $L_0'$-module. (The assertion is also
true for $p=5.$ \cite{BGP})
    When $p=3,$ the situation is more complicated.  In characteristic three,
there are series of simple graded Lie algebras which satisfy the
conditions of Kac's Recognition Theorem, but which are neither
classical Lie algebras nor Lie algebras of Cartan type. (See 
\cite{B}, \cite{Sk1}, \cite{St}.)  Moreover, for $q=1$, examples
exist in which $L_{-1}$ is not a restricted $L_0'-$module. All
simple depth-one graded Lie algebras over algebraically closed fields of characteristic three with
non-restricted $L_0'$-module $L_{-1}$ were determined in
\cite{BKK}. (In \cite{KO}, the authors classified all simple depth-one graded Lie algebras 
over algebraically closed fields of characteristic three in which $L_{-1}$
 is a {\it{restricted}} $L_0'-$module.) In \cite{BGK}, two-graded (i.e., depth-two, graded)
Lie algebras were examined, and it was proved that when $p = 3$
and $q = 2,$ the $L_0'-$module $L_{-1}$ must be restricted. For
$q=3,$ the corresponding statement was proved in \cite{GK1}.  It
was conjectured in \cite{BGK} that a non-restricted $L_0'-$module
$L_{-1}$ can exist only in Lie algebras of depth one. The present
paper completes the proof of that conjecture for $r > 1.$  We require that $r$ be greater than one in order to exclude, for example, $H(2:\mathbf{n}, \omega)$ with the reverse gradation; see also Example 0.3 below.  We note, as in \cite{BGK}, that because there are only
finitely many irreducible restricted modules for the derived
algebra of a classical reductive Lie algebra, what needs to be
considered in classifying graded Lie algebras over algebraically
closed fields of characteristic three is reduced. 
 In this paper, we prove the following theorem, which we will henceforth refer to as the ``Main Theorem."
\bigskip
\begin{Thm} {\bf (Main Theorem)} Let $L=L_{-q}\oplus L_{-q+1
}\oplus \cdots \oplus L_{-1} \oplus L_0\oplus L_1 \oplus \cdots
\oplus L_r, \, q
> 1, \, r > 1,$ be a finite-dimensional graded Lie algebra over an
algebraically closed field ${\mathbf F}$ of characteristic $p=3$
such that
\begin{enumerate}
\item[{\rm (A)}] $L_0$ is classical reductive;
\item[{\rm (B)}] $L_{-1}$ is an irreducible $L_0$-module (i.e.,
$L$ is {\it irreducible});
\item[{\rm (C)}] for all $j\geqq 0$, if $x\in L_j$ and
$[x,L_{-1}]=(0)$, then $x=0$ (i.e., $L$ is {\it transitive});

\item[{\rm (D)}] $L_{-i}=[L_{-i+1},L_{-1}]$ for all $i>1;$ and

\item[{\rm (E)}]$L_{-2} \nsubseteq M(L),$ where $M(L)$ is the
largest ideal of $L$ contained in the sum of the negative
gradation spaces. (See  Theorem \ref{Thm:1.3} below.)
\end{enumerate}
\noindent Then $L_{-1}$ is a restricted module for $L_{0}'$ under
the adjoint action of $L$ on itself.
\end{Thm}

To help to motivate hypothesis (E) above, we offer the following

\begin{Ex}  For characteristic $p \geqq 3,$ consider the irreducible transitive graded Lie
algebra
\begin{equation}R \eqdef {\mathcal O}(2: (1, \, 1)) \oplus H(2: (1, \, 1)) =
\bigoplus_{i = -2 - 2(p-1)}^{2p-5}R_i,\nonumber\end{equation}
\noindent where $R_i = H(2:(1, \, 1))_i$ for $i \geqq -1,$ and $R_i
= {\mathcal O}(2: (1, \, 1))_{i+2p}$ for $-2p = -2-2(p-1) \leqq i
\leqq -2.$ Here, the divided-power algebra ${\mathcal O}(2: (1, \,
1))$ is an abelian ideal of $R,$ and $H(2:(1, \, 1))$ has its
usual Lie algebra bracket operation and its usual action on ${\mathcal O}(2:
(1, \, 1)),$ except that $[D_{x_1}, \, D_{x_2}] =
x_1^{p-1}x_2^{p-1} \in R_{-2}.$ (See  Theorem \ref{Thm:1.3} below.) Then $R/M(R) \cong R/{\mathcal
O}(2: (1, \, 1))$ has depth
one. In general, if we consider the free Lie algebra generated by
the local part of any depth-one graded Lie algebra $L,$ and take a
co-finite-dimensional subideal $C$ of the maximal ideal $D$ in the
negative part (See  \cite{BW}.), then $M(L \oplus D/C) = D/C,$ and
$(L \oplus D/C)/M(L \oplus D/C) \cong L$ has depth one.

\end{Ex}
 To further illustrate the necessity of the requirement in the Main Theorem that $r$ be greater than one, we offer the following 
\begin{Ex}  Consider a graded Lie algebra
\begin{equation} L = \bigoplus_{i= - q}^1 L_i \nonumber\end{equation}
\noindent where
\begin{equation} L_1 = \langle \partial_j, \, j = 1, \dots n \rangle \nonumber \end{equation}
\begin{equation} L_0 = S + T \nonumber \end{equation}
\noindent where $S$ and $T$ are classical, and

\begin{equation} L_{-i} = S \otimes \langle x_1^{k_1} \dots x_n^{k_n}, \, k_1 + \dots \, + k_n = i \rangle, \nonumber\end{equation}
\noindent and where $T$ is a subalgebra of $\mathfrak{gl}(n) = \langle x_i\partial_j, \, i,j = 1, \dots, n\rangle = W(n:1)_0$ 
with a non-restricted action on $\langle x_1, \dots, x_n\rangle;$
i.e., the representation of $T$ in $\langle x_1, \dots, x_n\rangle$ is a non-restricted representation of $T.$ 
 (See  (ii) of Theorem \ref{Thm:1.3}.)  More succinctly, $L = S\otimes A(n:1) + 1\otimes T + 1\otimes \langle\partial_i, \, i = 1, \dots, n\rangle.$  
The minimal example corresponds to $T = \mathfrak{sl}(2)$ and $n = 3,$ so that, therefore, $q = 6.$
\end{Ex}

    {\it We have noted that the Main Theorem has been proved for $q=2$ in \cite{BGK}
and for $q=3$ in \cite{GK1}.  When we refer to the Main Theorem to
substantiate certain claims below, it will be for the cases
already proved.}

\bigskip

	We conclude this section with a sketch of the plan of the rest of the work.  In Section 1, 
we establish terminology and notation, and gather some previously known results.
We continue to gather previous results at the beginning of Section 2; here, however, the
results are of a more technical nature, and we immediately apply them 
to showing that under quite natural assumptions,
 the $L_0'-$module $L_{-2}$ is irreducible.  
We conclude Section 2 by establishing other technical lemmas 
that we'll use later on.  In Sections 3, 4, and 5, 
we prove the Main Theorem under the ``natural'' assumptions just referred to, 
for the cases $q \geqq 6,$ $q = 4,$ and $q = 5,$ respectively.  
In Section 6, we complete the proof of the Main Theorem in its full generality.

\setcounter{section}{0}
\begin{section}{Preliminaries} \end{section}

	In this section, we recall definitions and introduce notation, after which we gather
results from the literature that we will use later in the work.

    Recall that if one takes a $\mathbb Z$-form (Chevalley basis) of a complex simple Lie algebra and
reduces the scalars modulo $p,$ one obtains a Lie algebra over $I/(p).$  If ${\mathbf
F}$ is any field of characteristic $p$ ,  then, by tensoring the Lie algebra we obtained over $I/(p)$ by  ${\mathbf
F},$ we obtain a Lie algebra over  ${\mathbf
F}.$  In characteristic $p,$ any such Lie algebra so obtained is referred to as {\it{classical}}, 
even those with root systems E$_6$, E$_7$, E$_8$, F$_4$, and G$_2.$ This process may result in a Lie algebra with a non-zero
center; such a Lie algebra is still referred to as classical,
as is the quotient of such a Lie algebra by its center. For
example, the Lie algebras $\mathfrak{gl}(pk)$ and
$\mathfrak{pgl}(pk)$ are both considered to be classical Lie
algebras. Thus, a classical Lie algebra ${\mathfrak g}$ may have a
nontrivial center $\mathfrak{z}(\mathfrak{g}),$ as do the Lie algebras
$\mathfrak{gl}(pk), \mathfrak{sl}(pk),$ and, if $p = 3,$ E$_6.$ It
could also happen that a classical Lie algebra has a noncentral
ideal, as do the Lie algebras $\mathfrak{gl}(pk),$
$\mathfrak{pgl}(pk),$ and, if $p = 3,$ G$_2.$ In characteristic
three, G$_2$ contains an ideal $I$ isomorphic to
$\mathfrak{psl}(3),$ and $G_2/I \cong \mathfrak{psl}(3),$ as well.
 A {\it classical reductive Lie algebra\/} $\mathfrak{g}$ is the
sum of commuting ideals $\mathfrak{g}_j$ which are classical Lie
algebras, and an at-most-one-dimensional center
$\mathfrak{z}(\mathfrak{g})$:

\begin{equation}\mathfrak{g} = \mathfrak{g}_1 + \cdots + \mathfrak{g}_k + \mathfrak{z}(\mathfrak{g}) \label{eq:1.1}\end{equation}

\noindent For any classical Lie algebra $\mathfrak{g}_j,$ the
derived algebra $\mathfrak{g}_j'$ is a restricted Lie algebra with a natural $p$-structure such
that $e_{\alpha }^{[p]}$ $=0,$ and $h_{i}^{[p]}$ $=h_i$ for any
Chevalley basis $\{ e_{\alpha },\; h_{i}\; | \; \alpha$ $\in R,\;
i = 1,\dots, \text {rank}(\mathfrak{g}_j') \}$ of
$\mathfrak{g}_j'$, where $R$ is the root system of the
corresponding complex simple Lie algebra. For a classical
reductive Lie algebra ${\mathfrak g},$ we will consider only that
$p-$structure on ${\mathfrak g}'$ $=$ ${\mathfrak g}_1'$ $+$ $\cdots$
$+$ ${\mathfrak g}_k'$ whose restriction to each classical summand is the natural $p$-structure on that summand.

    Let $\pi:  L \longrightarrow {\mathfrak gl}(n)$ be a finite-dimensional irreducible
representation of a restricted Lie algebra $L$.  The {\it
character} of $\pi$ is the linear functional  $\chi$ on $L$ such that
$\chi (y)^pI = \pi (y)^p-\pi (y^{[p]})$ for all $y\in L$.   The
representation  $\pi$ is restricted when the character  $\chi$ equals $0$.

\bigskip

\begin{Lem}  \label{Lem:1.2} (See  \cite[Lemma 1]{BG}.)
Assume that $L$ is a graded Lie algebra satisfying conditions
(A)-(D) of the Main Theorem. If $\chi $ is the character of $L_0'$
on $L_{-1}$, then $L_0'$ has character $-j\chi $ on $L_j$ for all
$j$. \end{Lem}

 The following theorem of Weisfeiler \cite{W} plays a fundamental
r\^ole in the study of graded Lie algebras.  In what follows, we will sometimes
refer to Theorem \ref{Thm:1.3} as ``Weisfeiler's Theorem.''

\begin{Thm} {\bf (Weisfeiler's Theorem)} \label{Thm:1.3}  
Let $L=L_{-q}\oplus \cdots \oplus L_{-1} \oplus
 L_0 \oplus L_1 \oplus \cdots \oplus L_r$ be a graded Lie algebra such
that conditions (B)-(D) of the Main Theorem hold. Let $M(L)$
denote the largest ideal of $L$ contained in $L_{-q}\oplus \cdots
\oplus L_{-1}.$ Then

\begin{enumerate}
\item[{\rm (i)}] $L/M(L)$ is semisimple and contains a unique minimal ideal
$I=S\otimes {\mathcal O}(n:{\mathbf 1})$, where $S$ is a simple
Lie algebra, $n$ is a non-negative integer, and ${\mathcal
O}(n:{\mathbf 1})={\mathbf F}[x_1,\ldots ,x_n]/(x_1^p,\ldots
,x_n^p).$  The ideal $I$ is graded and $I_i=(L/M(L))_i$ for all
$i<0.$
\item[{\rm (ii)}] {\bf Degenerate Case} If $I_1=(0)$, then for some $\mathfrak{k}$, $1 \leqq \mathfrak{k}  \leqq
n$, the algebra ${\mathcal O}(n:{\mathbf 1})$ is graded by setting
$\deg (x_i)=-1$ for $1\leqq i\leqq \mathfrak{k}$  and $\deg (x_i)=0$ for
$\mathfrak{k}  < i\leqq n$.  Then $I_i=S\otimes {\mathcal O}(n:{\mathbf
1})_i$ for all $i$, \; $L_2=(0),$ \; $I_0=[L_{-1},L_1],$ and
$L_1\subseteq \{D\in 1\otimes \Der {\mathcal O}(n:{\mathbf 1})
\vert \deg (D)=1\}.$
\item[{\rm (iii)}] {\bf Non-Degenerate Case} If $I_1\ne (0)$, then $S$ is graded and $I_i=S_i\otimes
{\mathcal O}(n:{\mathbf 1})$ for all $i$. Moreover, $(0) \ne
[L_{-1},L_1] \subseteq I_0.$ 
\end{enumerate}
 \end{Thm}
Set
\begin{equation}L_{<0} \eqdef \bigoplus_{i = -q}^{-1}L_i,\nonumber\end{equation}
\noindent and
\begin{equation}L_{>0} \eqdef \bigoplus_{i = 1}^{r}L_i.\nonumber\end{equation}

\noindent  Set $L_{\leqq} \eqdef L_{<} + L_0,$ and $L_{\geqq} \eqdef L_{0} + L_{>}.$  

\begin{Thm} \label{Thm:1.4}  (See Theorem 0.1  of
\cite{GK2}.)  Let $L=\oplus _{i\in {\mathbb Z}} L_i$ be a
non-degenerate graded Lie algebra over an algebraically closed
field of characteristic $p>2$ satisfying conditions (A)-(D) of
the Main Theorem. If $[[L_{-1}, V], V] = 0$ for some proper $L_0$-submodule
$V \subset L_1$, $\dim V > 1,$ then $\dim L = \infty $. \end{Thm}
 
 By Lemma \ref{Lem:1.2}, the representation of $L_0'$ on $L_1$
 is restricted
 when and only when the representation of $L_0'$ on $L_{-1}$ is restricted.  
Since no non-restricted representation of $L_0'$ can have dimension one,
 we have the following corollary.

\begin{Cor}  \label{Cor:1.41}
Let $L$ be as in the above theorem, and suppose that the representation of $L_0'$ on $L_{-1}$ is not restricted.  If $[[L_{-1}, V], V] = 0$ for some proper $L_0$-submodule
$V \subset L_1$,  then $L$ is an infinite-dimensional Lie algebra.  \end{Cor}

    We will make use of the following results from \cite{BGK}.  For
 definitions of the Lie algebras $L(\epsilon),$ $M,$ $H(2:{\mathbf n},\omega),$ and $CH(2:{\mathbf n},\omega)$ mentioned in the
 conclusion of Proposition \ref{Pro:1.4} below, see, for example, Section 2 of \cite{BGK}.
 When we make use of certain properties of these Lie algebras in
 later sections, we will explicitly state the properties we need.

\begin{Pro}  \label{Pro:1.4} (See Lemma 2.12 of
\cite{BGK}.) Let $L=L_{-1}\oplus L_0\oplus L_1 \oplus \cdots
\oplus L_r$ be a graded Lie algebra satisfying conditions (A),
(B), and (C) of the Main Theorem, and suppose that $L_1 \neq 0.$
If $L_{-1}$ is a nonrestricted $L_0'-$module, then either $L$ is
isomorphic to one of the Lie algebras L$(\epsilon)$ or $M$, or $L$
is a Hamiltonian Lie algebra such that $H(2:{\mathbf n},\omega)
\subseteq L \subseteq CH(2:{\mathbf n},\omega),$ where ${\mathbf
n} = (1, n_2)$, $\omega = (\exp x^{(3)}) dx \wedge dy,$ and the
grading is of type (0,1).
\end{Pro}

\begin{Cor}  \label{Cor:1.5} (See Corollary 2.13 of
\cite{BGK}.)  Under the assumptions of Proposition \ref{Pro:1.4},
$L_0'$ $\cong \mathfrak{sl}(2),$ $L_1$ is an irreducible
three-dimensional $L_0'-$module, and $[L_1,$ $L_1]$ $= 0.$  In
addition, $[L_{-1},$ $L_1] \cong \mathfrak{sl}(2)$ if and only if $L$
is a Hamiltonian Lie algebra; otherwise, $[L_{-1}, \, L_1] \cong {\mathfrak{gl}}(2).$
\end{Cor}

\begin{Lem}  \label{Lem:1.6} (See Lemma 2.14 of
\cite{BGK}.) Let $L=L_{-1}\oplus L_0\oplus L_1 \oplus \cdots
\oplus L_q$ be one of the Lie algebras $L(\epsilon)$, $M,$ or
$H(2:{\mathbf n},\omega)$ with ${\mathbf n} = (1, n_2)$,   let $
\chi $ be the nonzero character of the $L_0$-module $L_{-1}$, and
let $V$ be an $L$-module such that $l^3 \cdot V=(0)$ for any $l\in
L_{-1}\cup L_1.$  Suppose that $W$ is an irreducible
$L_0$-submodule of $V$ with character $\chi_W = \zeta \chi, \, \zeta
\in {\mathbf{F}}^\times$, and suppose that $L_1\cdot W=(0).$  Then
$L_{-1}^2\cdot W\ne (0).$ Similarly, if $L_{-1}\cdot W=(0),$  then
$L_{1}^2\cdot W\ne (0).$
\end{Lem}

{\it In what follows, all Lie algebras will be finite-dimensional
over an algebraically closed field ${\mathbf F}$ of characteristic
$p=3$.  The commutator ideal $[L,L]$ of a Lie algebra $L$ will be
denoted by $ L',$ and the $i^{th}$ commutator $(\ad X)^{i-1}X$ of
any set $X$ will be written as $X^i.$  The annihilator of an
$L_0-$module $M \subseteq L$ in an $L_0-$module $N \subseteq L$
will be denoted by $\Ann_NM.$}

\begin{section}{Properties of irreducible transitive graded Lie algebras} \end{section}

	This section contains technical lemmas and a proof that under hypotheses which we list here, the
representation of the null component on the minus-two component must be irreducible.

    We begin this section by recalling a few results from \cite{BG}.
Let $L,$ $M(L),$ $I,$ and $S = \sum_{i = -q}^sS_i$ be as in Weisfeiler"s Theorem
(Theorem \ref{Thm:1.3}).  Throughout this section, we make the
following two blanket assumptions: 

\bigskip

\begin{enumerate}
\item[{\rm (i)}] $M(L) = 0$
\item[{\rm (ii)}] $I = S.$
\end{enumerate}

	In this regard, please see  \cite[(2.4.6)]{Sk2}, and note that the Lie algebra $S + L_0$ satisfies hypothesis (B) of the Main Theorem.

\begin{Lem}  \label{Lem:2.1} (See \cite[Lemma 6]{BG}.) For
any $x$ in $L \backslash L_{-q}, \, [L_{-1}, \, x] \neq 0.$
\end{Lem}

\begin{Lem}  \label{Lem:2.2} (See  \cite[Lemma 7]{BG}.) $S_{j} = (\ad
S_{-1})^{s-j}S_s$ for all $j, \, -q \leqq j \leqq s.$  If $q(t-1)
\leqq s,$ then $(\ad S_{-q})^tS = 0$ if and only if $(\ad
S_{-q})^tS_i = 0$ for some $i, \, q(t-1) \leqq i \leqq s.$
\end{Lem}

\begin{Lem}  \label{Lem:2.3} (See \cite[Lemma 8]{BG}.) $[L_{-q}, \, L_i] \neq
0$ for all $i = 0, \dots r.$  In addition, $S_j = (\ad
L_{-1})^{r-j}L_r$ for all $j, \, -q \leqq j \leqq r - 1,$ so that $s
= r - 1$ or $r.$ $S_s$ is an irreducible $S_0$-module.
\end{Lem}

\begin{Lem}  \label{Lem:2.4} (See \cite[Lemma 9]{BG}.) $S_{-q}$ is an
irreducible $S_0$-module. In particular, $L_{-q}$ is an
irreducible $L_0$-module.\end{Lem}

\begin{Lem}  \label{Lem:2.5} (See \cite[Lemma 10]{BG}.) $\Ann_{L_0}L_i
\cap \Ann_{L_0}V_{i+1} = 0$ for all $i = -q, \dots, r-1,$
where $V_{i+1}$ is any non-zero $L_0-$submodule of $L_{i+1}.$
\end{Lem}

\begin{Lem}  \label{Lem:2.6} (See \cite[Lemma 11]{BG}.) $\Ann_{L_i}L_{-q}
\cap \Ann_{L_i}V_{-q+1} = 0$ for all $i = 0, \dots, r,$
where $V_{-q+1}$ is any non-zero $L_0-$submodule of $L_{-q+1}.$
\end{Lem}

\begin{Lem}  \label{Lem:2.7} (See \cite[Lemma 12]{BG}.)
$\Ann_{L_{q-1}}L_{-q+1} = 0.$ \end{Lem}

\begin{Lem}  \label{Lem:2.8} (See \cite[Lemma 13]{BG}.) If $r \geqq q,$ then
$L_{-q+i} = [L_{-q}, \, L_i]$ for $i = 0, 1, \dots q-1.$
\end{Lem}

\begin{Lem}  \label{Lem:2.9} (See \cite[Lemma 14]{BG}.) Let $U$ and $V$ be
$L_0-$submodules of $L$ such that $[U, \, V] \subseteq L_0$ and
$[U, \, [U, \, V]] = 0.$  Then $\{\ad[u, \, v] | u \in U, \, v \in
V \}$ is weakly closed (in the sense of \cite[p.31]{J});
consequently, if $(\ad [u, \, v])^i M = 0, \, u \in U, \, v \in V,
\,$ for some $i > 1,$ and $L_0-$module $M,$ then $\ad_M[U, \, V]$
is ``associative nilpotent.'' (See  Theorem II.2.1 of \cite{J}.)
\end{Lem}

\begin{Lem}  \label{Lem:2.10} (See \cite[Lemma 15]{BG}.)   If $r \geqq q,$ then
$(\ad L_{-q})^2  L  \ne (0).$
 \end{Lem}

\begin{Cor}  \label{Cor:2.11}  If $s \geqq q,$ then $[S_{-q}, \, [S_{-q},
\, S_q]]  \ne (0).$   In particular, if $s \geqq q,$ then $[L_{-q},
\, [L_{-q}, \, L_q]]  \ne (0).$ \end{Cor}

{\bf Proof}  By Lemma \ref{Lem:2.10}
applied to the Lie algebra $S+L_0,$ $(\ad S_{-q})^2 S\neq 0.$  Thus, in view of Lemma \ref{Lem:2.2}, it cannot be that $(\ad S_{-q})^2 S_q = 0.$  $~~\square$

\begin{Lem}  \label{Lem:2.12} (See  \cite[Lemma 16]{BG}.) Let $V$ be an
$L_0-$submodule of $L_{-q+i}$ for some $i,$ where $0$ $<$ $i$ $\leqq$ $\frac{q}{2},$ and
suppose $[V, \, L_{q-i-1}] = 0 = [V, \, [V, L_{q-i}]].$
Suppose further that $L_{-q+i-1}$ is an irreducible
$L_0-$submodule of $L,$ and that $[L_{-q+i-1}, \, L_{q-i}]$ $\neq$ $0$
(so that it equals $L_{-1}).$  Then $V$ $=$ $0.$
\end{Lem}

\begin{Lem}  \label{Lem:2.13} (See \cite[Lemma 17]{BG}.) Suppose that $V$ is an
irreducible \linebreak $L_0-$submodule of $L_{-q+i}$ for some $i,$ where $0$ $< i$ $<
\frac{q-1}{2},$ such that $[V,$ $L_{q-i-1}]$ $\neq 0$ (so that it
equals $L_{-1}$).  Then $L_{-q+i}$ is an irreducible $L_0-$module;
i.e., it equals $V.$
\end{Lem}

From our observations at the beginning of this section, we have
that $S\subseteq L\subseteq \Der S$ where $S=S_{-q}\oplus S_{-q+1}
\oplus \cdots \oplus S_{-1}\oplus S_0 \oplus S_1\oplus \cdots
\oplus S_s$ is a simple Lie algebra with $S_i=L_i$ for $i<0.$
Since $S_1$ is an $L_0$-submodule of $L_1,$ it follows that if
$L_1$ is an irreducible $L_0-$module, then $S_1=L_1.$ If, in
addition, $L$ is generated by its local part $L_{-1} \oplus L_0
\oplus L_1,$ then for $i \geqq 1,$ we have
\begin{equation}L_i = L_1^i = S_1^i \subseteq S_i \subseteq
L_i,\nonumber\end{equation}
\noindent so that $S_i=L_i$ for $i>0,$ and $L$ could differ from
$S$ only in the null component. In particular, $s$ would equal
$r.$  (See  Lemma \ref{Lem:2.3}.)

    In the lemmas that follow, we will consider graded Lie algebras
$L= L_{-q}\oplus L_{-q+1}\oplus \cdots \oplus L_{-1}\oplus
L_0\oplus L_1\oplus \cdots \oplus L_r$ satisfying assumptions (i)
and (ii) below. Other assumptions will be noted in the statements
of the results for which we use them. Note, for example, that, as
noted in the paragraph above, assumption (iii) follows from
assumptions (iv) and (v) (and assumptions (i) and (ii), of
course). Also, assumption (viii) can be assumed whenever the
previous assumptions are true, since if they hold, we can reverse
the gradation, and have that all of the hypotheses of the Main
Theorem continue to be true for the reversed gradation.  Indeed, (iv) 
is the ``reverse'' of hypothesis (B) of the Main Theorem, (vi) is the ``reverse'' of (C), 
and (v) is the ``reverse'' of (D).  Of course, by the transitivity (C) of $ L$,
 there can be no ideals of $L$ in the positive part of $L$, so (E) holds 
in the ``reverse'' direction, also  In addition, by Lemma \ref{Lem:1.2},
 the representation of $L_0'$ on L$_{-1}$ is restricted
 if and only if the representation of $L_0'$ on $L_1$ is restricted, so (vii) holds in the ``reverse'' direction, as well.
\begin{enumerate}
\item[{\rm (i)}] $L$ satisfies conditions (A) - (E) of the Main
Theorem.
\item[{\rm (ii)}] $L\subseteq \Der S$ where $S=S_{-q}\oplus
S_{-q+1}\oplus S_{-1}\oplus S_0\oplus S_1\oplus \cdots \oplus S_s$
is a simple graded Lie algebra.
\item[{\rm (iii)}] $L_i=S_i,\ i\ne 0.$
\item[{\rm (iv)}]  $L_1=S_1$  is an irreducible $L_0$-module. (See  the discussion before \eqref{eq:6.01}.)
\item[{\rm (v)}] $L_{i+1}=[L_i,L_1]$ for $i>0$.
\item[{\rm (vi)}] If $x$ is a non-zero element in $L_{-i}$ for
some $i \geqq 0$,  then $[L_1,x] \ne (0).$
\item[{\rm (vii)}] The character $\chi $ of $L_0'$ on $L_{-1}$ is
non-zero.
\item[{\rm (viii)}] $r\geqq q.$
\end{enumerate}
 \begin{Lem}  \label{Lem:2.14}     If assumptions (iv)  and (v) hold, and $S_1$ $\neq 0,$ then $\Ann_{L_0}L_1$ $= 0.$
\end{Lem}

{\bf Proof} Suppose, on the contrary, that $A_0 \eqdef
\Ann_{L_0}L_1 \neq 0.$  Then (as in \cite[Lemma 18]{BG}) we have by transitivity (C) and irreducibility (B) that
\begin{equation}[L_{-1}, \, L_1]  = [[L_{-1}, \, A_0], \, L_1] = [[L_{-1}, \, L_1], \, A_0] \subset A_0,\nonumber\end{equation}
\noindent so that by transitivity (C),
\begin{equation}0 \neq [S_{-1}, \, L_1] \subseteq A_0 \cap S_0 \subseteq  \Ann_{S_0}S_1.\nonumber\end{equation}
\noindent  Since for $i > 0,$ we have by assumptions (iv)  and (v)
that $S_i = S_1^i,$ and since $J \eqdef S_{-q} \oplus S_{-q+1}
\oplus \cdots \oplus S_{-1} + A_0 \cap S_0$ is invariant under
$\ad{S_i}, \, -q \leqq i \leqq s,$ it follows, from our assumption that $S_1 \neq 0,$ 
that $J$ is a proper ideal of the simple Lie algebra $S;$ i.e., we have obtained a
contradiction. Thus, we must conclude that $\Ann_{L_0}L_1 =
0.$ $~~\square$
\begin{Lem}  \label{Lem:2.15} If assumption (vi) holds, then $[V_{-2},$
$L_1]$ $= L_{-1}$ for any non-zero $L_0-$submodule $V_{-2}$ of
$L_{-2}.$
\end{Lem}

{\bf Proof} This lemma follows from assumptions (vi) and (B).
$~~\square$

\begin{Lem}  \label{Lem:2.16} If assumption (iv)  holds, then
$\Ann_{L_i}L_{-q} = 0$ for all $i > 0.$
\end{Lem}

{\bf Proof}  Consider first the case in which $i = 1.$  If
$\Ann_{L_1}L_{-q} \neq 0,$ then, since we are assuming (iv) 
that $L_1$ is an irreducible $L_0$-module, we would have
$\Ann_{L_1}L_{-q} = L_1.$  But then

\begin{equation}[L_{-q}, \, L_1] = [L_{-q}, \, \Ann_{L_1}L_{-q}] = 0,\nonumber\end{equation}

\noindent to contradict Lemma 2.3.  Consequently,
$\Ann_{L_1}L_{-q}$ $= 0.$  Now, if $Q_i$  \linebreak  $\eqdef
\Ann_{L_i}L_{-q}$ $\neq 0$ for some $i$ $> 1,$ then by
transitivity (C), we would have
\begin{equation}0 \neq (\ad L_{-1})^{i-1}Q_i \subset \Ann_{L_1}L_{-q},\nonumber\end{equation}
\noindent to contradict what we just showed.  Thus,
$\Ann_{L_i}L_{-q} = 0$ for all $i > 0,$ which is what we
wanted to show. $~~\square$

\begin{Lem}  \label{Lem:2.17}  If assumptions (vi) and (viii) hold, then
$L_{-q}$ $= [L_{-q+i},$ $L_{-i}],$ $0$ $\leqq i$ $\leqq q.$
\end{Lem}

{\bf Proof}  By Lemmas \ref{Lem:2.3} and \ref{Lem:2.4}, $[L_0, \,
L_{-q}] = L_{-q},$ so the lemma is true for $i = 0$ and $i = q.$
For $i = 1,$ we use Lemmas \ref{Lem:2.1} and \ref{Lem:2.4}.    Now
note that for $1 \leqq i \leqq q-1,$ we have
\begin{eqnarray} [(\ad{L_1})^{i}L_{-q}, \, (\ad{L_1})^{q-i}L_{-q}]
&=&
[(\ad{L_1})^{i}L_{-q}, \, [L_1, \, (\ad{L_1})^{q-(i+1)}L_{-q}]\nonumber\\
&=& [(\ad{L_1})^{i+1}L_{-q}, \,
(\ad{L_1})^{q-(i+1)}L_{-q}],\nonumber
\end{eqnarray}
\noindent so that in view of (vi), (B), and  Lemmas \ref{Lem:2.1} and \ref{Lem:2.4},
\begin{eqnarray} L_{-q} &=& [(\ad{L_1})^{1}L_{-q}, \, (\ad{L_1})^{q-1}L_{-q}]\nonumber\\
&=& [(\ad{L_1})^{i}L_{-q}, \, (\ad{L_1})^{q-i}L_{-q}], \, 1 \leqq i
\leqq q.\nonumber
\end{eqnarray}
\noindent  Then
\begin{equation}L_{-q} = [(\ad{L_1})^{i}L_{-q}, \, (\ad{L_1})^{q-i}L_{-q}]
\subseteq [L_{-q+i}, \, L_{-i}] \subseteq L_{-q} , \, 1 \leqq i
\leqq q.\nonumber\end{equation}
\noindent $\square$

\begin{Lem}  \label{Lem:2.18} Let $V$ be any (non-zero)
irreducible $L_0-$submodule of $L_{-q+i},$ where  $0$ $\leqq i$ $< \frac{q-1}{2}.$  If  assumptions (v) and
(viii) hold, then $[V,$
$L_{q-(i+1)}]$ $\neq 0,$ so that $[V,$ $ L_{q-(i+1)}]$ $= L_{-1}.$
Moreover, $L_{-q+i}$ is an irreducible $L_0-$module.
\end{Lem}

{\bf Proof}  If $[V, \, L_{q-(i+1)}] = 0,$ then, since the
positive gradation spaces are assumed (v) to be generated by
$L_1,$

\begin{eqnarray} [V, \, L_{q-i}] &=& [V, \, [L_1, \, L_{q-(i+1)}]]\nonumber\\
&=& [[V, \, L_1], \, L_{q-(i+1)}]\nonumber\\
&\subseteq& [L_{-q+(i+1)}, \, L_{q-(i+1)}].\nonumber
\end{eqnarray}

\noindent Consequently, since $i$ is assumed to be less than
$\frac{q-1}{2},$ so that $2i+1 < q,$

\begin{eqnarray} [V, \, [V, \, L_{q-i}]] &\subseteq&
[V, \, [L_{-q+(i+1)}, \, L_{q-(i+1)}]]\nonumber\\
&=& [L_{-q+(i+1)}, \, [V, \,L_{q-(i+1)}]]\nonumber\\
&=& [L_{-q+(i+1)}, \, 0] = 0.\nonumber
\end{eqnarray}

\noindent Now, $L_{-q} = S_{-q}$ is an irreducible $L_0-$module by Lemma
\ref{Lem:2.4}, so we can assume by induction on $i$ that $L_{-q+i-1}$ is
an irreducible $L_0-$module.  Also, since $[L_{-q}, \, L_{q-1}] =
L_{-1}$ by Lemma \ref{Lem:2.3} and (B), we can assume by induction
that $[L_{-q+i-1}, \, L_{q-i}]$ $=$  $L_{-1}.$  But then Lemma
\ref{Lem:2.12} would imply that $V$ $= 0,$ contrary to assumption.
Thus, $[V,$ $ L_{q-(i+1)}]$ is a non-zero $L_0-$submodule of
$L_{-1},$ so that by irreducibility (B), $[V,$ $L_{q-(i+1)}]$ $=
L_{-1}.$ The last assertion follows from Lemma \ref{Lem:2.13}.
$~~\square$
\begin{Lem}  \label{Lem:2.19}  Suppose that assumptions (vi) and (viii)
hold.  Then for any $i, 1 \leqq i \leqq q,$ we have  and $[L_{-q},
\, [L_{-i}, \, L_i]] = L_{-q};$ in particular, $[L_{-i}, \, L_i]
\neq 0.$
\end{Lem}

{\bf Proof}  The lemma will follow from Lemma \ref{Lem:2.4}  once
we show that $[L_{-q},$  $[L_{-i}, \, L_i]] \neq 0.$   For $i =
q,$ the lemma follows from Corollary \ref{Cor:2.11}. Let $1 \leqq i
< q.$ Then by Lemma \ref{Lem:2.17}, we have $L_{-q} = [L_{-q+i},
\, L_{-i}],$ and by Lemma \ref{Lem:2.8}, we have $L_{-q+i} =
[L_{-q}, \, L_{i}].$ Then we have

\begin{equation}L_{-q} = [L_{-q+i}, \, L_{-i}] = [[L_{-q}, \, L_{i}], \,L_{-i}]
= [L_{-q}, [L_{-i}, \,  L_{i}]],\nonumber\end{equation}

\noindent so that $[L_{-q}, [L_{-i}, \,  L_{i}]] \neq 0,$ as
required. $~~\square$

\begin{Lem}  \label{Lem:2.20}  Suppose that assumptions (iv) , (v), (and
therefore, (iii)), (vi), and (viii) hold and that $0 < i < \frac{q
- 3}{2}.$ Then $[L_{-q+i}, \ L_{q-i+1}] = L_{1}.$
\end{Lem}

{\bf Proof} Suppose that $[L_{-q+i}, \ L_{q-i+1}] = 0.$  Then,
since by (iii) and Lemma \ref{Lem:2.3}, $L_j = S_j = [S_{j+1}, \,
S_{-1}] = [L_{j+1}, \, L_{-1}] \subseteq L_j$ for all $j,\,  0 < j < r,$
we have

\begin{eqnarray}
[L_{-q+i}, \, L_{q-i}] &=& [L_{-q+i}, \, [L_{q-i+1}, \,
L_{-1}]]\nonumber\\
&=& [L_{q-i+1}, \, [L_{-q+i}, \, L_{-1}]]\nonumber\\
&=& [L_{q-i+1}, \, L_{-q+i-1}],\nonumber
\end{eqnarray}

\noindent so that (since $i < \frac{q-3}{2}$ implies that $2i + 3
< q,$ so that {\it a fortiori} $2i - 1 < q)$
\begin{eqnarray}
[L_{-q+i}, \ [L_{-q+i}, \, L_{q-i}]] &=& [L_{-q+i}, \ [L_{q-i+1},
\,
L_{-q+i-1}]]\nonumber\\
&=& [L_{-q+i-1}, \ [L_{-q+i}, \, L_{q-i+1}]]\nonumber\\
&=& [L_{-q+i-1}, \ 0] = 0.\nonumber
\end{eqnarray}

\noindent Let $v \in L_{-q+i} $ and $u \in L_{q-i}.$  Then (since $2i + 1 < 2i + 3 < q,$ so that $[v, \, L_{-q+i+1}] = 0$),
\begin{eqnarray} 2(\ad [v, \, u])^2L_{-q+i+1} &=& (\ad v)^2(\ad
u)^2L_{-q+i+1}\nonumber\\
&\subseteq& (\ad v)^2L_{q-i+1}\nonumber\\
&\subseteq& (\ad v)[L_{-q+i}, \ L_{q-i+1}] = (\ad v)\cdot 0 =
0.\nonumber
\end{eqnarray}

\noindent Consequently, $\ad_{L_{-q+i+1}}[L_{-q+i}, \, L_{q-i}]$
is a nilpotent set of linear transformations by Lemma
\ref{Lem:2.9}.  Since we are assuming that $i < \frac{q-3}{2},$ we
have $i+1 < \frac{q-1}{2},$ so we can apply Lemma \ref{Lem:2.18}
to conclude that $L_{-q+i+1}$ is an irreducible $L_0-$module. It
follows that $\ad_{L_{-q+i+1}}[L_{-q+i}, \, L_{q-i}]$ annihilates
$L_{-q+i+1}.$ Thus (since, again,  $2i + 1 < 2i + 3 < q)$

\begin{equation}0 = [[L_{-q+i}, \, L_{q-i}], \, L_{-q+i+1}] = [L_{-q+i}, \, [L_{q-i}, \,
L_{-q+i+1}]].\nonumber\end{equation}
\noindent If $[L_{q-i}, \, L_{-q+i+1}]$ $\neq 0,$ then, since
$L_1$ is assumed (iv)  to be irreducible, $[L_{q-i},$ $L_{-q+i+1}]
$ would have to equal $L_1,$ and the above-displayed formula would
imply a lack of $\{ 1 \}-$transitivity (vi) of $L$ in its negative
part. It follows that $[L_{q-i}, \, L_{-q+i+1}]$ $= 0.$  Then, in
view of our initial assumption that $[L_{-q+i}, \ L_{q-i+1}]$ $=$ $0,$
we would have

\begin{equation}0 = [[L_{-q+i}, \, L_{1}], \, L_{q-i}] = [[L_{-q+i}, \, L_{q-i}], \,
L_{1}],\nonumber\end{equation}

\noindent to contradict Lemma \ref{Lem:2.14}, in view of Lemma
\ref{Lem:2.19}. Thus, it must be that $[L_{-q+i},$ $L_{q-i+1}]$
$\neq 0,$ so that by the assumed irreducibility (iv) of $L_1,$
$[L_{-q+i}, \ L_{q-i+1}]$ $= L_{1},$ as required.  $~~\square$

\begin{Lem}  \label{Lem:2.21}  Let $q > 5,$ and suppose that assumptions
(iv) , (v) (and therefore (iii)), (vi), and (viii) hold. If $q$ is
even, then $L_{-2}^{\frac{q}{2}} = L_{-q},$ while if $q$ is odd,
then $L_{-2}^{\frac{q-1}{2}} = L_{-q+1}.$
\end{Lem}

{\bf Proof} We have by Lemma \ref{Lem:2.8}  that $L_{-2} =
[L_{-q}, \, L_{q-2}]$ and by Lemma \ref{Lem:2.18} (since $q > 3)$
that $L_{-1} = [L_{-q+1}, \, L_{q-2}].$ Thus, for any $j, \, 1 < j
< q-1,$ and any non-zero $L_0-$submodule $V_{-j}$ of
$L_{-j},$ we have by $\{ -1 \}-$transitivity (Lemma \ref{Lem:2.1})
that

\begin{equation}0 \neq [V_{-j}, \, L_{-1}] = [V_{-j}, \, [L_{-q+1}, \,
L_{q-2}]] = [L_{-q+1}, \, [V_{-j}, \,
L_{q-2}]].\nonumber\end{equation}

\noindent  Consequently, $[V_{-j},$ $L_{q-2}]$ $\neq 0.$  Then by
Lemma \ref{Lem:2.16} when $j$ $< q - 2,$ or, when $j$ $= q - 2,$
by Lemma \ref{Lem:2.19},

\begin{eqnarray}0 \neq [L_{-q}, \, [V_{-j}, \,  L_{q-2}]]
&=& [V_{-j}, \, [L_{-q}, \, L_{q-2}]]\nonumber\\
&=& [V_{-j}, \, L_{-2}] \nonumber
\end{eqnarray}
\noindent by Lemma \ref{Lem:2.8}.  If, for $j = 1, \, 2, \, \dots, \, [\frac{q}{2}] - 1,$ we successively let $V_{-2j}
\eqdef L_{-2}^j,$ we can conclude that $L_{-2}^{\frac{q}{2}} \neq
0$ if $q$ is even, and
$L_{-2}^{\frac{q-1}{2}} \neq 0$ if $q$ is odd. Then, by Lemmas \ref{Lem:2.4} and
\ref{Lem:2.18}, respectively, $L_{-2}^{\frac{q}{2}} = L_{-q}$ if $q$ is even, or 
$L_{-2}^{\frac{q-1}{2}} = L_{-q+1}$ if $q$ is odd. $~~\square$

\begin{Lem}  \label{Lem:2.22}  Let $q > 5,$ and suppose that assumptions
(iv) , (v) (and therefore (iii)), (vi) and (viii) hold.  Then
$L_{-2}$ is an irreducible $L_0-$module.
\end{Lem}

{\bf Proof}       Let $V_{-2}$ be any irreducible $L_0-$submodule
of $L_{-2}.$ Since $[L_{-q+1},$ $[V_{-2},$ $L_{q}]]$ $= [V_{-2},$
$[L_{-q+1},$ $L_{q}]]$ $= [V_{-2},$ $L_1]$ $= L_{-1}$ by Lemmas
\ref{Lem:2.20} and \ref{Lem:2.15}, it follows that for any $j,$ $0
< j < \frac{q}{2}$ (i.e., $0 < j$ $\leqq \frac{q-1}{2}$) for which
$V_{-2}^j$ $\neq 0,$ we have by transitivity (Lemma \ref{Lem:2.1})
that

\begin{equation}[L_{-q+1}, \, [V_{-2}^j, \, [V_{-2}, \, L_{q}]]] = [V_{-2}^j,
[V_{-2}, \, [L_{-q+1}, \, L_{q}]]= [V_{-2}^j, L_{-1}] \neq 0,\nonumber\end{equation}

\noindent so we conclude that $[V_{-2}^j, \, [V_{-2}, \, L_{q}]]$
$\neq 0$ and therefore that  $(\ad V_{-2})^{j+1} L_{q} \neq 0.$   Thus, so long
as  $2(j+1) < q$ (i.e., $j < \frac{q}{2} - 1),$ we have by Lemma
\ref{Lem:2.16} that
\begin{eqnarray} 0 &\neq& [L_{-q}, \, [V_{-2}^j, \, [V_{-2}, \, L_{q}]]]\nonumber\\
&=& [V_{-2}^j, \, [V_{-2}, \, [L_{-q}, \, L_{q}]]]\nonumber\\
&\subseteq& V_{-2}^{j+1}.\nonumber
\end{eqnarray}

\noindent Thus, $V_{-2}^j \neq 0$ for all $j,$ $0 < j \leqq
\frac{q-1}{2},$ and $(\ad V_{-2})^jL_q \neq 0$ for all $j,$ $0 < j
\leqq \frac{q+1}{2}.$  If $q$ is odd, then, since $q
> 5,$ we have by Lemma \ref{Lem:2.18} that $V_{-2}^{\frac{q-1}{2}}
= L_{-q+1},$ while if $q$ is even, we have $V_{-2}^{\frac{q}{2}-1}
= L_{-q+2}.$

    In the case of odd $q,$  we have, by the irreducibility (B) of $L$, that
$L_{-1}$  $= (\ad V_{-2})^\frac{q+1}{2}L_q,$ so that
\begin{eqnarray} L_{-2} &=& [L_{-1}, \, L_{-1}]\nonumber\\
&=& [L_{-1}, \, (\ad V_{-2})^\frac{q+1}{2}L_q]\nonumber\\
&\subseteq&[V_{-2}, \, [L_{-q}, \, L_q]] + [V_{-2}, \, L_0]\nonumber\\
&\subseteq& V_{-2}.\nonumber
\end{eqnarray}

\noindent Thus, when $q$ is odd, we see that $L_{-2}$ is
irreducible.
    In the case of even $q,$ we have by Lemma \ref{Lem:2.18} (since $q > 5)$ that $L_{-1}$
$= [L_{-q+2}, \, L_{q-3}]$ $= [V_{-2}^{\frac{q}{2}-1}, \,
L_{q-3}]$ $\subseteq (\ad V_{-2})^{\frac{q}{2}-1}L_{q-3}$
$\subseteq L_{-1}.$  By (D) and transitivity (Lemma
\ref{Lem:2.1}),
\begin{eqnarray} L_{-q+1} &=& [L_{-1}, \, L_{-q+2}]\nonumber\\
&=& [(\ad V_{-2})^{\frac{q}{2}-1}L_{q-3}, \, L_{-q+2}]\nonumber\\
&=& [(\ad V_{-2})^{\frac{q}{2}-1}L_{q-3}, \, V_{-2}^{\frac{q}{2}-1}]\nonumber\\
&\subseteq& (\ad V_{-2})^{q-2}L_{q-3}\nonumber\\
&\subseteq& L_{-q+1},\nonumber
\end{eqnarray}
\noindent so that (See  also Lemma \ref{Lem:2.18}.) $(\ad
V_{-2})^{q-2}L_{q-3}$ $= L_{-q+1}.$ Then, by Lemma \ref{Lem:2.17},
\begin{eqnarray} L_{-q} &=& [L_{-1}, \, L_{-q+1}]\nonumber\\
&=& [L_{-1}, \, (\ad V_{-2})^{q-2}L_{q-3}]\nonumber\\
&\subseteq& [V_{-2}, \, L_{-q+2}] + [(\ad V_{-2})^{q-2}L_{-1}, \, L_{q-3}]\nonumber\\
&\subseteq& [V_{-2}, \, V_{-2}^{\frac{q}{2}-1}] + [0, \, L_{q-3}]\nonumber\\
&\subseteq& V_{-2}^{\frac{q}{2}}
\nonumber
\end{eqnarray}

\noindent Now, by Lemma \ref{Lem:2.16} and irreducibility (B), we
have
\begin{equation}L_{-1} = [L_{-q}, \, L_{q-1}] = [V_{-2}^{\frac{q}{2}}, \,
L_{q-1}] \subseteq (\ad V_{-2})^{\frac{q}{2}}L_{q-1} \subseteq
L_{-1}.\nonumber\end{equation}
\noindent Consequently, we have
\begin{eqnarray} L_{-2} &=& [L_{-1}, \, L_{-1}]\nonumber\\
&=& [L_{-1}, \, (\ad V_{-2})^{\frac{q}{2}}L_{q-1}]\nonumber\\
&=& [(\ad V_{-2})^{\frac{q}{2}}L_{-1},\, L_{q-1}] + [V_{-2}, \, L_{0}]\nonumber\\
&\subseteq& [0, \, L_{q-1}] + V_{-2}\nonumber\\
&\subseteq& V_{-2}\nonumber
\end{eqnarray}
\noindent as required.  $~~\square$

        Now suppose that (vii) holds.   Note that by (D) of the Main Theorem, $L_{-2} = [L_{-1}, \, L_{-1}];$ that is, $L_{-2}$ is spanned by brackets of elements of the three-dimensional $L_0-$module $L_{-1}.$ Consequently, $L_{-2}$ is at most three-dimensional.    On the other hand, by Lemma \ref{Lem:1.2}, the character $\chi$ is non-zero on $L_{-2},$ so $L_{-2}$ is not a restricted $L_0-$module, so its dimension is, in fact, three, as are all irreducible non-restricted ${\mathfrak{sl}} (2)-$modules in characteristic three.  Thus (Compare Lemma \ref{Lem:2.22}.), we have

\begin{Lem} \label{Lem:6.001} If assumption (vii) holds, then $L_{-2}$ is an irreducible three-dimensional non-restricted $L_0'-$module. \end{Lem}

\begin{Lem}  \label{Lem:2.23}  If $q > 2$ and assumptions  (vi) and (vii) hold, then 
$\Ann_{L_1}L_{-2}$ $= 0.$
\end{Lem}

{\bf Proof}  Set $A_1 = \Ann_{L_1}L_{-2},$ and suppose that
$A_1 \neq 0.$  Since
\begin{equation}[L_{-2}, \, [L_{-q+1}, \, A_1]] = [L_{-q+1}, \, [L_{-2}, \,
A_1]]= 0,\nonumber\end{equation}
\noindent we have
\begin{eqnarray} 0 &=& [L_{-2}, \, [L_{-q+1}, \, A_1]]\nonumber\\
&\supseteq& [[L_{-3}, \, L_1], \, [L_{-q+1}, \, A_1]] = [L_{-3}, \, [L_1, \, [L_{-q+1}, \, A_1]]]\nonumber\\
&\supseteq& [[L_{-4}, \, L_1], \, [L_1, \, [L_{-q+1}, \, A_1]]] = [L_{-4}, \, [L_1, \, [L_1, \, [L_{-q+1}, \, A_1]]]]\nonumber\\
&\cdots\nonumber\\
&\supseteq& [L_{-q+1}, \, (\ad L_1)^{q-3}[L_{-q+1}, \,
A_1]].\nonumber
\end{eqnarray}
\noindent Now, if $[L_{-q+1}, \, A_1] \neq 0,$ then by (vi) and
irreducibility (B), we would have $(\ad L_1)^{q-3}$ $[L_{-q+1}, \,
A_1]$ $= L_{-1},$ so that $[L_{-q+1}, \, L_{-1}] = 0,$ to
contradict transitivity (Lemma \ref{Lem:2.1}).  Thus, we must have

\begin{equation}\label{eq:2.23.1}[L_{-q+1}, \, A_1] = 0.\end{equation}

    Now, since $[L_{-q},$ $[A_1,$ $A_1]]$ $\subseteq [L_{-q+1},$
$A_1]$ $= 0,$ and, clearly, $[L_{-q+1},$ $[A_1,$ $A_1]]$ $= 0,$ it
follows from Lemma \ref{Lem:2.6} that $[A_1,$ $A_1]$ $= 0.$ Then
${L}^{\dagger}$ $\eqdef$ $(L_{-q}$ $\oplus$ $\cdots$ $\oplus$ $L_{-1}$
$\oplus$ $L_0$ $\oplus$ $A_1)/(L_{-q}$ $\oplus$ $\cdots$ $\oplus$ $L_{-2})$
 is a depth-one Lie algebra which satisfies conditions (A)
through (C) of the Main Theorem. Consequently, by Proposition
\ref{Pro:1.4}, $({L}^{\dagger})'$ is one of the Lie algebras
enumerated in the hypothesis of Lemma \ref{Lem:1.6}. If we set $V$
$=$ $L_{-q}$ $\oplus$ $L_{-q+1},$ then the $({L}^{\dagger})'-$module $V$
satisfies the hypotheses of Lemma \ref{Lem:1.6}.  Set $W$ $=
L_{-q+1},$ and note that by \eqref{eq:2.23.1}, $[A_1,$ $W]$ $=$ $0.$ 

	 If $(\zeta=) \, q - 1\not\equiv 0$ (mod 3), then Lemma \ref{Lem:1.6} would imply that  $[L_{-1},$ $[L_{-1},$ $L_{-q+1}]]$ $\neq$ $0.$ Since this is false, we must have $\zeta = q-1$ $\equiv 0$ mod 3.  On the other hand, if we then
set $W$ $=$ $L_{-q},$ we have $[L_{-1}, \, L_{-q}] = 0$ and, by \eqref{eq:2.23.1},

\begin{equation} [A_1, \,  [A_1, \, L_{-q}]] \subseteq [L_{-q+1}, \, A_1] = 0\nonumber\end{equation}

\noindent so we must, by similar reasoning,  conclude that $(\zeta=) \, q \equiv 0$ mod 3. Since both
$q-1$ and $q$ cannot be equivalent to zero modulo three, we have
arrived at a contradiction. We therefore conclude that
$\Ann_{L_1}L_{-2}$ $= A_1$ $= 0,$ as required. $~~\square$

\begin{Lem}  \label{Lem:2.26} If conditions (iv) , (vi)  and (vii) hold, then  
\begin{equation} \Ann_{L_1}L_{-2} = 0.\nonumber\end{equation}
\end{Lem}

{\bf Proof}  If $q = 2$, this lemma follows from Lemma \ref{Lem:2.16}.    If $q > 2,$ it follows from Lemma \ref{Lem:2.23}.
 \qed

\begin{Lem}  \label{Lem:2.28} If $M_1$ is a non-zero $L_0-$submodule of $L_1$ such that
$\Ann_{L_0}M_1$ $\neq 0,$ then $[[L_{-1}, M_1], \, M_1]$ $= 0,$ and $[M_1, \, M_1]$ $= 0.$
\end{Lem}

{\bf Proof}  Set $X \eqdef \Ann_{L_0}M_1 \neq 0.$  Then by
transitivity (C) and irreducibility (B), $[L_{-1}, \, X] = L_{-1},$ so
\begin{equation} [L_{-1}, \, M_1] = [[L_{-1}, \, X], \, M_1] = [[L_{-1}, \, M_1], X] \subseteq X. \nonumber \end{equation}
\noindent Thus, $[L_{-1}, [M_1, \, M_1]] \subseteq [[L_{-1},\,  M_1], \, M_1] \subseteq [X, \, M_1] = 0,$ so, by transitivity (C),
$[M_1, \, M_1] = 0.$ \qed

\begin{Lem}  \label{Lem:2.31} We may assume that  $(S_s =)\, [[L_{-1}, \, L_1], \, L_r] = L_r.$ \end{Lem}

{\bf Proof} Suppose $[[L_{-1}, \, L_1], \, L_r] \neq L_r.$  We distinguish two cases:

\begin{enumerate}
\item[\rm{(i)}] $[[L_{-1}, \, L_1], \, L_r] = 0.$
\item[\rm{(ii)}] $[[L_{-1}, \, L_1], \, L_r] \neq 0 .$
\end{enumerate}

	(i) Suppose first that $[[L_{-1}, \, L_1], \, L_r]$ were equal to zero.  Then we would have 

\begin{equation}  0 =[[L_{-1}, \, L_1], \, L_r]= [ L_1,[L_{-1} \, L_r]] \nonumber\end{equation}

\noindent  so $\sum_{i\geqq 0}(\ad L_{-1})^i[L_{-1}, \, L_r] \subseteq S$ would be an ideal in $S$ entailing equality by the definition of $S;$ in particular, we would have $[L_{-1}, \, L_r] = S_s.$ 

 If, in addition,  $[[L_{-1}, \, L_1], \, [L_{-1}, \, L_r]]$ were also equal to zero, we could repeat the argument and get that $\sum_{i\geqq 0}(\ad L_{-1})^i[L_{-1}, \, [L_{-1}, \, L_r]$ would be a proper ideal of $S$, to contradict the simplicity of $S$.  
We conclude that  $[[L_{-1}, \, L_1],$ $[L_{-1}, \, L_r]] \neq 0,$ and, in the case that  $[[L_{-1}, \, L_1], \, L_r] = 0$, that $[L_{-1}, \, L_r] = S_s,$ which is
an irreducible $L_0-$module by Lemma \ref{Lem:2.3}.

Then,

\begin{equation}[[L_{-1}, L_1], \, [L_{r}, \, L_{-1}]]  =  [[L_{-1}, L_1], \, S_s]] = S_s =[L_{-1}, L_r] \nonumber\end{equation} 

\noindent Hence, if in the case that  $[[L_{-1}, \, L_1], \, L_r] = 0$, we replace $L_r$ with $[L_{-1}, L_r]$, the lemma follows. 

	(ii) Now suppose that  $[[L_{-1}, \, L_1], \, L_r] \neq 0 .$  Since

\begin{equation} 0 \neq  [[L_{-1}, \, L_1], \, L_r]  =  [[L_{-1}, \, L_r], \, L_1] \subseteq [S_{r-1}, \, L_1] \subseteq S, \nonumber\end{equation}

\noindent it would follow that $ \sum_{i\geqq 0}(\ad L_{-1})^i[[L_{-1}, \, L_1],  \, L_r]$ would be an ideal of $S$ and hence all of $S.$ 
 In particular, $s$ would equal $r,$ and $[[L_{-1}, \, L_1],  \, L_r]$ would equal $S_s.$  Thus,  $[[L_{-1}, \, L_1], \, L_r] = S_s$ is an irreducible $L_0-$module by Lemma \ref{Lem:2.3}.
Consequently, if we replace $L$ by the Lie algebra generated by $L_{-1},$ $L_0,$ $L_1,$  and  $[[L_{-1}, \, L_1], \, L_r]$ $ (= [S_{s-1}, L_1] = S_s),$
then the highest gradation space will be of the form  $[[L_{-1}, \, L_1], \, L_r],$ as required.  \qed

\begin{Lem}  \label{Lem:2.25}  Let $L$ be as in the statement of the Main Theorem, and
suppose that $L_2$ $\neq 0,$ that $[L_{-2}, \, L_1]$ $= 0$ $=
[L_{-2},$ $\, L_2],$ and that assumption (vii) holds. Let
$\tilde{L}$ be the Lie subalgebra of $L$ generated by $L_{-1},$
$L_0,$ and $L_1.$  If $M(\tilde{L})$ is as in Weisfeiler's Theorem (Theorem
\ref{Thm:1.3}), then $\tilde{L}/M(\tilde{L})$ is Hamiltonian, and
we have $[L_{-1}, \, L_1] \cong \mathfrak{sl}(2).$
\end{Lem}

{\bf Proof} Let $\tilde{\tilde{L}}$ be the Lie subalgebra of $L$
generated by $L_{-1},$ $L_0,$ $L_1,$ and $L_2.$  Since $[L_{-2},$
$\, L_1]$ $= 0$ $= [L_{-2},$ $\, L_2],$ we have $M(\tilde{\tilde{L}})$ $=
L_{-q}$ $\oplus$ $\cdots$ $\oplus$ $L_{-2}$ $=
M(\tilde{L}).$ But the depth
$\tilde{\tilde{L}}/M(\tilde{\tilde{L}})$ is then one, so by
Proposition \ref{Pro:1.4},
$\tilde{\tilde{L}}/M(\tilde{\tilde{L}})$ is either Hamiltonian
(i.e., between $H(2:\underline{n}, \, \omega)$ and
$CH(2:\underline{n}, \, \omega))$ or is isomorphic to a Lie
algebra of type $L(\epsilon)$ or $M.$  However, the height of the
latter two Lie algebras is one, and, since $L_2 \neq 0,$ the
height of $\tilde{\tilde{L}}/M(\tilde{\tilde{L}})$ is at least
two.  Thus, $\tilde{\tilde{L}}/M(\tilde{\tilde{L}})$ must be
Hamiltonian, so $[L_{-1}, \, L_1] \cong \mathfrak{sl}(2).$  It now
follows from Corollary \ref{Cor:1.5} that $\tilde{L}/M(\tilde{L})$
is Hamiltonian, as well.  $~~\square$

\begin{Lem}  \label{Lem:2.255}  Let $L$ be as in the statement of the Main Theorem, 
and suppose that $\tilde{L}/M(\tilde{L})$ (as above) is isomorphic to $L(\epsilon)$ or $M$.  
Then $\Ann_{L_2}L_{-2} = 0$ and in the proof of the Main Theorem, 
where we assume that $[L_{-2}, \, L_1] = 0,$ 
we may assume that $L_2$ is an irreducible $L_0-$module.
\end{Lem}

{\bf{Proof}}  Suppose first that $\Ann_{L_2}L_{-2} \neq 0.$  
Then, as in the proof of the previous lemma, 
we may consider the Lie algebra $\tilde{\tilde{\tilde{L}}}$ 
generated by $L_{-1},$ $L_0,$ $L_1,$ and $\Ann_{L_2}L_{-2}.$  
Then $\tilde{\tilde{\tilde{L}}}/M(\tilde{\tilde{\tilde{L}}})$ has depth one 
and height greater than one, but null component isomorphic 
to ${\mathfrak{gl}}(2)$, to contradict Corollary \ref{Cor:1.5}.  
Thus, $\Ann_{L_2}L_{-2} = 0$, so that if $M_2$ 
is any irreducible $L_0-$module of $L_2,$ then $[L_{-2}, \, M_2] \neq 0.$  
Replacing $L$ by the Lie algebra generated by $L_{-1},$ $L_0,$ $L_1,$ and $M_2,$ 
we complete the proof of the lemma. $~~\square$

\begin{section}{The Main Theorem under additional assumptions}\end{section}

    In this and the following two sections, we assume that assumptions ((i) and (ii), of course),
(iv) , (v) (and therefore (iii)), (vi), and (viii) of the previous
section hold, so that, in particular, by Lemma \ref{Lem:2.22},
$L_{-2}$ is an irreducible $L_0-$module.
    We begin this section by forming the irreducible, transitive Lie algebra $B(L_{-2}).$
(See, for example, Section 3 of \cite{BG}.)  Indeed, consider the
subalgebra 
\begin{equation}E = E_{-\lfloor \frac q2  \rfloor} \oplus \cdots
\oplus E_0 \oplus \cdots \oplus E_{\lfloor \frac r2  \rfloor}\nonumber\end{equation} 
\noindent of
$L$ consisting of the gradation spaces $E_i = L_{-2}^{-i}$ for $i
< 0,$ and $E_i = L_{2i}$ for $i \geqq 0.$ Set $T_0= $
$\Ann_{E_0}E_{-1}= $ $\Ann_{L_0}L_{-2},$ and for
$i=1,2, \dots$,  let
\begin{equation} T_i=\{x\in E_i  \mid  [x,E_{-1}]\subseteq  T_{i-1}\}.\nonumber\end{equation}
\noindent  Then   
\begin{equation} {\mathcal T}=T_0\oplus T_1\oplus \ldots \oplus
T_{\lfloor \frac r2  \rfloor}\nonumber\end{equation}
\noindent   is an ideal of $E$,  and the
factor algebra   

\begin{equation}G= E /{\mathcal T} = G_{-{\lfloor \frac q2
\rfloor}}\oplus \cdots \oplus G_{-1} \oplus G_0\oplus G_1\oplus
\cdots \oplus G_{\lfloor \frac r2 \rfloor}\nonumber\end{equation}
\noindent is a transitive graded
Lie algebra (See  \cite[Lemma 3]{BG}.) Thus, the Lie algebra
$B(L_{-2}) \eqdef G$ satisfies conditions (A) - (D) of the Main
Theorem.  (It is shown in, for example, \cite{BGP} that the
process of forming $B(L_{-2})$ preserves condition
 (A).)

\bigskip

\begin{Lem}  \label{Lem:3.1}  Let $L$ be as in the
statement of the Main Theorem, and suppose that assumptions (i) -
 (viii) hold. If $q \geqq 6,$ then $B(L_{-2})$ is an
irreducible, transitive graded Lie algebra of height greater than or equal to 
two and depth greater than two, and $B(L_{-2})_{-2} \not \subseteq
M(B(L_{-2})).$ Consequently, since the depth of $B(L_{-2})$ is no
greater than half of the depth of $L,$ we can, using induction,
apply the Main Theorem to conclude that the character of the
representation of $B(L_{-2})_0'$ on $B(L_{-2})_{-1}$ is equal to
zero. Then the character of $L_0' = B(L_{-2})_{0}'$ on $L_{-1}$ is also zero.
\end{Lem}

{\bf Proof}    By Lemma  \ref{Lem:2.21}, when $q$ is
even, $L_{-2}^{\frac{q}{2}} = L_{-q},$ while when $q$ is odd,  $L_{-2}^{\frac{q-1}{2}} = L_{-q+1},$ so, in either case, the depth of $B(L_{-2})$ is greater than or equal to three.   In the even case, since by Corollary \ref{Cor:2.11} $[L_{-q}, \, [L_{-q}, \, L_q]]  \ne (0),$ we have

\begin{equation}\label{eq:3.01} 0 \neq [L_{-q}, \, [L_{-q}, \, L_q]] = [L_{-2}^{\frac{q}{2}}, \, [L_{-2}^{\frac{q}{2}}, \, L_q]] \subseteq (\ad L_{-2})^qL_q\end{equation}

\noindent so the height of  $B(L_{-2})$ is greater than or equal to three, since $r \geq q \geq 6.$ 

	In the odd case, since $q \geqq 7,$ we have by Lemma \ref{Lem:2.18} that $[L_{-q+2},$ $L_{q-3}]$ $= L_{-1}.$  If $[L_{-q+1},$ $L_{q-3}]$ were equal to zero, then we would have

\begin{eqnarray}
 0 &=& [L_{-q+2}, \, 0]\nonumber\\
&=& [L_{-q+2}, \, [L_{-q+1}, \, L_{q-3}]]\nonumber\\
&=& [L_{-q+1}, \, [L_{-q+2}, \, L_{q-3}]] = [L_{-q+1},
\,L_{-1}]\nonumber
\end{eqnarray}

\noindent to contradict transitivity (Lemma \ref{Lem:2.1}).  We
therefore conclude that 

\begin{equation}\label{eq:3.02}[L_{-q+1}, \, L_{q-3}] \neq 0.\end{equation}

\noindent Then $0 \neq [L_{-q+1}, \, L_{q-3}] = [L_{-2}^{\frac{q-1}{2}}, \, L_{q-3}]\subseteq (\ad L_{-2})^{\frac{q-1}{2}}L_{q-3}$ so the height of $B(L_{-2})$ is at least two in the odd case.

      We must now verify hypothesis (E) of the Main Theorem for
the Lie algebra $B(L_{-2});$ that is, we must show that $[L_{-2},$
$L_{-2}]$ is not contained in $M(B(L_{-2})).$  Thus, suppose that
$[L_{-2},$ $L_{-2}]$ is contained in $M(B(L_{-2}));$ i.e., that 

\begin{equation}\label{eq:3.07}[L_2, \, [L_{-2}, \,L_{-2}]] = 0\end{equation}

\noindent We will arrive at a contradiction by successively considering the two
cases:

\begin{enumerate} 
\item[{\rm 1)}] even $q$, and
\item[{\rm 2)}]  odd $q.$
\end{enumerate}

\begin{enumerate}
\item[{\rm 1)}] Suppose first that $q$ is even.  Then by Lemma \ref{Lem:2.21},
$L_{-q}$ $= L_{-2}^{\frac{q}{2}}$ \linebreak $=$ $(\ad
L_{-2})^{\frac{q}{2}-2}[L_{-2},$ $L_{-2}],$ so that $L_{-q}$ $\subseteq
M(B(L_{-2})).$  Then by Lemma \ref{Lem:2.8}, ($B(L_{-2})_{-1}$ $=)
L_{-2}$ $= [L_{-q},$ $L_{q-2}]$ $\subseteq M(B(L_{-2})),$ so that
we would have by the definition of $M(B(L_{-2}))$ that

\begin{equation}[L_{-2}, \, L_2] \subseteq L_0 \cap M(B(L_{-2})) = 0,\nonumber\end{equation}

\noindent to contradict, for example, \eqref{eq:3.01} above. Thus,
$q$ cannot be even.

\item[{\rm 2)}]  Next suppose that $q$ is odd. Then we have by Lemma \ref{Lem:2.21} again,
\begin{equation} L_{-q+1} = L_{-2}^{\frac{q-1}{2}} = (\ad
L_{-2})^{\frac{q-1}{2}-2}[L_{-2}, 
  L_{-2}],\nonumber\end{equation}

\noindent so that $L_{-q+1}$ $\subseteq M(B(L_{-2})).$  Since by Lemma \ref{Lem:2.22}, $L_{-2}$ is an irreducible $L_0-$module,
it follows from \eqref{eq:3.02} that $L_{-2}$ $= [L_{-q+1},$ $L_{q-3}]$ $\subseteq
M(B(L_{-2})).$ But we saw at the conclusion of 1) above that
$L_{-2}$ cannot be contained in $M(B(L_{-2})).$  This second
contradiction shows that $B(L_{-2})_{-2}$ $= [L_{-2},$ $L_{-2}]$
is in fact not contained in $M(B(L_{-2})),$ no matter what the
parity of $q$ is.
\end{enumerate}
    Consequently, we can conclude that $B(L_{-2})$ satisfies
hypothesis (E), and therefore all of the hypotheses of the Main
Theorem.  Since the depth of $B(L_{-2})$ is greater than one, 
but less than or equal to $\frac{q}{2},$  
we can now apply to conclude that the character
$\chi$ of $B(L_{-2})_0'$ on $B(L_{-2})_{-1}$ is zero.  Then
$\frac{1}{2}\chi,$ which is the character of $L_0'$ on $L_{-1},$
must be zero as well, and Lemma \ref{Lem:3.1} is proved.
$~~\square$

\bigskip

    We now address the depth-four and depth-five
cases individually.

\begin{section}{The depth-four case} \end{section}

    Suppose

\begin{equation}L = L_{-4} \oplus L_{-3} \oplus \cdots \oplus L_{r}\nonumber\end{equation}

\noindent satisfies conditions (i) through (viii) of Section 2.
    By Lemma \ref{Lem:2.15},

\begin{equation}[L_{-2}, L_1] = L_{-1}.\nonumber\end{equation}

    By Lemma \ref{Lem:2.17}, $[L_{-2}, L_{-2}]$ $= L_{-4}.$
Furthermore, from Lemma \ref{Lem:2.3} we have $[L_{-4},$ $L_{2}]$
$\neq 0.$ Then

\begin{equation}\label{eq:4.1} 0 \neq [L_{-4}, L_{2}] = [[L_{-2}, L_{-2}], L_{2}]
\subseteq [L_{-2}, \, [L_{-2}, L_{2}]]\end{equation}

    Now let $V_{-2}$ be any irreducible $L_0-$submodule of
$L_{-2}.$ If $[V_{-2},$ $L_3] = 0,$ then by Lemma \ref{Lem:2.3}
and irreducibility (B), $0$ $= [L_{-4},$ $[V_{-2},$ $L_3]]$ $=
[V_{-2},$ $[L_{-4},$ $L_3]]$ $= [V_{-2},$ $L_{-1}]$ to contradict
transitivity (Lemma \ref{Lem:2.1}).  Thus, we can assume that
$[V_{-2},$ $L_3]$ $= L_1,$ since we are assuming that $L_1$ is
$L_0-$irreducible (iv). Then by Lemma \ref{Lem:2.15}, $L_{-1}$ $=
[V_{-2},$ $ [V_{-2},$ $L_3]].$ Then we have by condition (D) of
the Main Theorem that

\begin{eqnarray} L_{-2} &=& [L_{-1}, \, L_{-1}]\nonumber\\
&=& [L_{-1}, \, [V_{-2}, \, [V_{-2}, \, L_3]]]\nonumber\\
&\subseteq& [[L_{-1}, \, V_{-2}], \, [V_{-2}, \, L_3]] + [V_{-2}, \, L_0]\nonumber\\
&=& [V_{-2}, \, [[L_{-1}, \, V_{-2}], \, L_3]] + [V_{-2}, \, L_0]\nonumber\\
&\subseteq&  [V_{-2}, \, L_0]\nonumber\\
&\subseteq&  V_{-2}\nonumber
\end{eqnarray}

\noindent so that $L_{-2}$ is an irreducible $L_0-$module.

    Thus, in the depth-two irreducible, transitive graded Lie
algebra $B(L_{-2}),$ we have by \eqref{eq:4.1} that $B(L_{-2})_1$
$\neq 0.$ Furthermore, it follows again from \eqref{eq:4.1} above
that
\begin{equation}B(L_{-2})_{-2} = L_{-4} = [L_{-2}, \, L_{-2}] \not \subseteq
M(B(L_{-2})),\nonumber\end{equation}
\noindent so that hypothesis (E) of the Main Theorem is satisfied
for $B(L_{-2}),$ as are the other hypotheses of the Main Theorem.
Then the Main Theorem (proved for the case $q = 2$ in \cite{BGK})
applies to show that the representation of $B(L_{-2})_0'$ on
$B(L_{-2})_{-1}$ is restricted. Consequently, the character $\chi$
of $L_0'$ on $L_{-2}$ is zero, as must be the character
$\frac{1}{2}\chi$ of $L_0'$ on $L_{-1}.$ 

\bigskip

\begin{section}{The depth-five case} \end{section}

    Suppose

\begin{equation}L = L_{-5} \oplus L_{-4} \oplus \cdots \oplus L_{r}\nonumber\end{equation}

\noindent satisfies conditions (i) through (viii) of Section 2.
    Since $L$ is transitive (C), $[L_{-1}, \,L_5] \neq 0.$  Hence, by Lemma \ref{Lem:2.16}, 
$[L_{-5}, \, [L_{-1}, \,L_5]] \neq 0,$ so that, by irreducibility (B), $[L_{-5},
\, [L_{-1}, \,L_5]] = L_{-1}.$

    Now suppose that $[L_{-4}, \, L_5] = 0.$  If also $[L_{-3}, \,
[L_{-3}, \, L_5]]$ $= 0,$ then we would have

\begin{equation}0 = [L_{-1}, \, [L_{-3}, \, [L_{-3}, \, L_5]]]
= [L_{-3}, \, [L_{-3}, \, [L_{-1}, \, L_5]]].\nonumber\end{equation}

    However, since $[L_{-5}, \, [L_{-1}, \, L_5]]$ $= L_{-1},$
we have by transitivity (Lemma \ref{Lem:2.1}) that

\begin{equation}0 \neq [L_{-1}, \, L_{-3}] = [L_{-3}, \, [L_{-5}, \, [L_{-1}, \, L_5]]]
= [L_{-5}, \, [L_{-3}, \, [L_{-1}, \, L_5]]]\nonumber\end{equation}

\noindent so that $[L_{-3}, \, [L_{-1}, \, L_5]]$ $\neq 0.$ Since
$L_1$ is assumed (iv)  to be irreducible, we must have $[L_{-3},$
$[L_{-1},$ $L_5]]$ $= L_1.$  Then by $\{ 1 \}-$transitivity (vi),

\begin{equation}0 \neq [L_{-3}, \, L_1] = [L_{-3}, \, [L_{-3}, \, [L_{-1}, \, L_5]]],\nonumber\end{equation}

\noindent contrary to what was derived above.  Thus, we can assume
that $[L_{-3},$ $[L_{-3},$ $L_5]]$ $\neq 0.$  But then, by the
irreducibility (B) of $L$, we must have $[L_{-3},$ $[L_{-3},$
$L_5]]$ $= L_{-1}.$  Then, by $\{ -1 \}-$transitivity (Lemma
\ref{Lem:2.1}),

\begin{equation}0 \neq [L_{-4}, \, L_{-1}]
= [L_{-4}, \, [L_{-3}, \, [L_{-3}, \, L_5]]] = [L_{-3}, \,
[L_{-3}, \, [L_{-4}, \, L_5]]],\nonumber\end{equation}

\noindent so that $[L_{-4}, \, L_5]$ $\neq 0.$
    Since we are assuming (iv)  that $L_1$ is irreducible, it follows that

\begin{equation}\label{eq:5.1}[L_{-4}, \, L_5] =
L_1.\end{equation}

    Now let $V_{-2}$ be any non-zero $L_0-$submodule of $L_{-2}.$    Then
we have by \eqref{eq:5.1} and $\{ 1 \}-$transitivity (vi) that

\begin{equation}0 \neq [V_{-2}, \, L_1] = [V_{-2}, \, [L_{-4}, \, L_5]]
= [[L_{-4}, \, [V_{-2}, \,  L_5]],\nonumber\end{equation}

\noindent so

\begin{equation}[L_{-4}, \, [V_{-2}, \,  L_5]] = L_{-1},\nonumber\end{equation}

\noindent by the irreducibility (B) of $L.$ Then, by transitivity (Lemma \ref{Lem:2.1}),

\begin{equation}0 \neq [V_{-2}, \,  L_{-1}] = [V_{-2}, \,  [L_{-4}, \, [V_{-2}, \,  L_5]]]
= [L_{-4}, \, [V_{-2}, \,  [V_{-2}, \,   L_5]]],\nonumber\end{equation}

\noindent so that $[V_{-2}, \,  [V_{-2}, \,   L_5]]$ $\neq 0.$
Thus, by the assumed irreducibility (iv)  of $L_1,$ we must have
$[V_{-2}, \,  [V_{-2}, \,   L_5]]$ $= L_1.$  Then, as above, by
the $\{ 1 \}-$transitivity (vi) of $L,$ we have

\begin{equation}0 \neq [V_{-2}, \, L_1] = [V_{-2}, \,  [V_{-2}, \,  [V_{-2}, \,   L_5]]],\nonumber\end{equation}

\noindent so that

\begin{equation}\label{eq:5.2}[V_{-2}, \,  [V_{-2}, \,  [V_{-2}, \,   L_5]]] =
L_{-1},\end{equation}

\noindent by the irreducibility (B) of $L.$ Then, since the
negative gradation spaces are generated (D) by $L_{-1},$ we have

\begin{equation}\label{eq:5.11}L_{-2} = [L_{-1}, \, L_{-1}] =
[L_{-1}, \, [V_{-2}, \,  [V_{-2}, \,  [V_{-2}, \,   L_5]]]]
\subseteq [V_{-2}, \,   L_0] \subseteq V_{-2}, \end{equation}

\noindent so that $L_{-2}$ is an irreducible $L_0-$module.

    Now, by Lemma \ref{Lem:2.17}, $L_{-5} = [L_{-2}, \, L_{-3}].$
Consequently, in view of (D) and \eqref{eq:5.2} above

\begin{eqnarray} L_{-5} &=& [L_{-2}, \, L_{-3}]\nonumber\\
&=& [L_{-2}, \, [L_{-1}, \, L_{-2}]]\nonumber\\
&=& [L_{-2}, \, [L_{-2}, \, [L_{-2}, \, [L_{-2}, \, [L_{-2}, \,
L_5]]].\nonumber \end{eqnarray}

\noindent Then by $\{ 1 \}-$transitivity (vi), we have

\begin{equation}0 \neq
[L_1, \, L_{-5}] = [L_1, \, [L_{-2}, \, [L_{-2}, \, [L_{-2}, \,
[L_{-2}, \, [L_{-2}, \, L_5]]] \subseteq [L_{-2}, \, [L_{-2}, \,
L_0]],\nonumber\end{equation}

\noindent so that $[L_{-2}, \, L_{-2}] \neq 0.$  Furthermore,
since the negative gradation spaces are generated (D) by $L_{-1},$
we have by \eqref{eq:5.2} above that

\begin{eqnarray} L_{-4} &=& [L_{-1}, L_{-3}]\nonumber\\
&=& [L_{-1}, [L_{-1}, L_{-2}]]\nonumber\\
&=& [L_{-1}, \, [L_{-2}, \, [L_{-2}, \, [L_{-2}, \, [L_{-2}, \,
L_5]]]]]\nonumber\\
&\subseteq& [L_{-2}, \, L_{-2}],\nonumber \end{eqnarray}

\noindent so

\begin{equation}\label{eq:5.12}L_{-4} = [L_{-2}, \, L_{-2}].\end{equation}

    Now suppose that $[L_{-4}, \, L_2]$ $= 0,$ and suppose further that $[L_{-3}, \, L_2]$
$= 0.$  Then we would have by (C) and (D) that
\begin{equation}0 = [L_{-4}, \, L_2] = [[L_{-3}, \, L_{-1}], \, L_2] =
[L_{-3}, \, [L_{-1}, \, L_2]] = [L_{-3}, \, L_1]\nonumber\end{equation}
\noindent by the assumed irreducibility (iv)  of $L_1,$ to
contradict $\{ 1 \}-$transitivity (vi). Thus, $[L_{-3}, \, L_2]$
$\neq 0,$ so that by the irreducibility (B) of $L,$ $[L_{-3}, \,
L_2]$ $= L_{-1}.$ Then

\begin{equation}[[L_{-3}, \, [L_{-4}, \, L_2]]
= [[L_{-4}, \, [L_{-3}, \, L_2]] = [L_{-4}, \, L_{-1}] \neq 0,\nonumber\end{equation}

\noindent by $\{-1 \}-$transitivity (Lemma \ref{Lem:2.1}).  Thus,
it must be true that

\begin{equation}\label{eq:5.205}[L_{-4}, \, L_2] \neq 0,\end{equation}

\noindent so that

\begin{equation}\label{eq:5.3}0 \neq [L_{-4}, \, L_2] = [[L_{-2}, \, L_{-2}], \, L_2]
\subseteq [L_{-2}, \, [L_{-2}, \, L_2]].\end{equation}

  	By \eqref{eq:5.11} and construction, $B(L_{-2})$ is a transitive, irreducible depth-two Lie algebra. By
\eqref{eq:5.3} above, $B(L_{-2})_1$ $\neq 0.$ By \eqref{eq:5.12}
and \eqref{eq:5.205} above, $B(L_{-2})$ satisfies hypothesis (E)
of the Main Theorem. Therefore, as in the depth-four case above,
it follows from the Main Theorem (proved for the case $q = 2$ in
\cite{BGK}) that the character of $B(L_{-2})_{0}'$ on
$B(L_{-2})_{-1}$ $= L_{-2}$ is zero. Consequently, the character
of $L_0'$ on $L_{-1}$ is zero, as well, and $L_{-1}$ is a
restricted $L_0'-$module. 

\bigskip

\begin{section}{Conclusion of the proof of the Main Theorem} \end{section}
    Let $L$ be as in the statement of the Main Theorem, and let $S = \sum_{i = -q}^sS_i$ be as in Weisfeiler's Theorem (Theorem \ref{Thm:1.3}).    If $[S_{>0}, L_{-2}]$ were equal to zero, then $L_{-q} \oplus \cdots \oplus  L_{-2} =  \sum_{i \geqq 0}(\ad (L_{-1}))^{i}L_{-2}$ would be an ideal of the simple Lie algebra $S$.  Consequently, it must be that $[S_{>0}, L_{-2}] \neq 0.$  Let $j > 0$ be minimal such that
 
\begin{equation}\label{eq:6.00001} [L_{-2}, \, S_j] \neq 0. \end{equation}

\noindent We wish to show that $j = 1.$  Suppose not.   

 	We begin by establishing a few basic properties. By the Jacobi Identity,
if $1 \leq k \leq j-1$ and $2 \leq i \leq k + 1,$ then (since in the sum that follows, $0 \leqq \kappa \leqq i - 2 \leqq k - 1 \leqq j - 2,$ so that $j - 1 \geqq k \geqq k - \kappa \geqq 1;$ here $\kappa$ is the number of $(\ad L_{-1})$s that act on $S_k$ before it brackets with $L_{-2}$ and annihilates it)

\begin{eqnarray} [L_{-i}, \, S_{k}] &=& [(\ad (L_{-1}))^{i-2}L_{-2}, \,
S_{k}]\nonumber\\
&\subseteq& \sum_{0\leq\kappa
\leq i-2}(\ad (L_{-1}))^{(i-2-\kappa)}[L_{-2}, \, S_{k-\kappa}]\nonumber\\
&=& 0;\nonumber
\end{eqnarray}
\noindent i.e.,
\begin{equation}\label{eq:6.1}[L_{-i}, \, S_{k}] = 0, \, 2 \leq i \leq k + 1, \, 1 \leq
k \leq j-1. \end{equation}
\noindent In particular, we have (letting $i = k$, $i = k+1$, and $k = j-1$, respectively)

\begin{equation}\label{eq:6.1.1}[L_{-k}, \, S_k] = 0, \, 2 \leqq k \leqq j-1,\end{equation}

\begin{equation}\label{eq:6.1.2}[L_{-k-1}, \, S_k] = 0, \, 1 \leq k \leq j-1, \hbox{ and }\end{equation}

\begin{equation}\label{eq:6.1.3}[L_{-i}, \, S_{j-1}] = 0, \, 2 \leq i \leq j.\end{equation}

We will now show that $[L_{-i}, \, S_j] \neq 0, \, 1 \leq i \leq j+1.$
Since $L$ is irreducible (B), it will follow that

\begin{equation}\label{eq:6.2}[L_{-j-1}, \, S_j] = L_{-1}.\end{equation}

\noindent Note that (since we are assuming that
$[L_{-2}, \, S_j] \neq 0$ and $[L_{-2}, \, S_{j'}]$ = $0$ for $1
\leq j' \leq j-1),$ we have by (D), by \eqref{eq:6.1.3}, by
transitivity (C), and by induction, that

$$[L_{-i}, \, S_{j}] = [[L_{-i+1}, \, L_{-1}], \, S_{j}] = [[L_{-i+1}, \,
S_{j}], \, L_{-1}] \neq 0, \, 2 < i \leq j+1.$$

\noindent Thus, we have (in view of transitivity (C) and \eqref{eq:6.00001}) that

\begin{equation}\label{eq:6.3.1}[L_{-i}, \, S_j] \neq 0, \, 1 \leq i \leq j+1. \end{equation}

	We will now show that

\begin{equation}\label{eq:6.3.11} \Ann_{L_j}L_{-j} =  \Ann_{L_j}L_{-2}\end{equation}

\noindent Thus, suppose first that $T_j$ is an $L_0-$submodule of $L_j$ such that $[L_{-j}, \, T_j] = 0.$  Then, in view of $D$ and \eqref{eq:6.1.3}, we have 

\begin{equation}\label{eq:6.3.111}0 = [L_{-j}, \, T_j] = [[L_{-j+1}, \, L_{-1}], \, T_j] = [[L_{-j+1}, \, T_j], \,  L_{-1}]\end{equation}

\noindent so that by transitivity (C), we have $[L_{-j+1}, \, T_j] = 0.$  Then we may replace $L_{-j}$ in \eqref{eq:6.3.111} by, successively, $L_{-j + 1},$  $L_{-j + 2},$ etc., to arrive at $[L_{-2}, \, T_j] = 0.$  On the other hand, if we rather define $T_j = \Ann_{L_j}L_{-2},$ then, again in view of  \eqref{eq:6.1.3}, we may bracket the equation $[L_{-2}, \, T_j] = 0$ with $L_{-1}$ again and again to conclude that

\begin{equation} [L_{-3}, \, T_j] = 0, \, [L_{-4}, \, T_j] = 0, \, \cdots, [L_{-j}, \, T_j] = 0,\nonumber\end{equation}

\noindent as required to establish \eqref{eq:6.3.11}.

\bigskip

\begin{Lem} \label{Lem:6.011} If $j$ is as above, then $[L_{-j-1}, \, S_{j+1}] \neq 0.$\end{Lem}

{\bf Proof} Suppose that $[L_{-j-1}, \, S_{j+1}] = 0$.  Then, since by \eqref{eq:6.2} $[L_{-j-1}, \, S_{j}]$ $=$ $L_{-1},$ we have

\begin{equation} [[S_{j}, \,  S_{j+1}], \, L_{-j-1}] = [[L_{-j-1}, \, S_{j}], \, S_{j+1}] = [L_{-1}, \, S_{j+1}] = S_j \neq 0.\nonumber\end{equation}

\noindent Consequently, $[S_{j}, \,  S_{j+1}] \neq 0.$ Then we have

\begin{eqnarray} [[[S_{j}, \,  S_{j+1}],  \,  S_{j+1}], \, L_{-j-1}] &=& [[[L_{-j-1}, \, S_{j}], \, S_{j+1}], \,  S_{j+1}]\cr &=& [[L_{-1}, \, S_{j+1}], \,  S_{j+1}] =  [S_j, \,  S_{j+1}] \neq 0.\nonumber\end{eqnarray}

\noindent Continuing in this manner, we get homogeneous non-zero $L_0-$submodules of arbitrarily high gradation degree.  This of course cannot happen in a finite-dimensional Lie algebra, so it must be that $[L_{-j-1}, \, S_{j+1}]\neq 0$. \qed

\begin{Lem}  \label{Lem:6.02} If  $j \geqq 2$, and $k$ is the smallest integer greater than one for which $[L_{-k}, \, S_r] \neq 0,$ then $k = j.$ \end{Lem}

{\bf Proof} Since by (D), $L_{<0}$ is generated by $L_{-1}$, we have by  \eqref{eq:6.3.1} and \eqref{eq:6.1.1} and (iv)  that

\begin{equation} (0 \neq) [L_{-j}, \, S_j] = [[L_{-1}, L_{-j+1}], \, S_j] = [L_{-1}, \,[ L_{-j+1}, \, S_j]] = [L_{-1}, \, L_1]\nonumber\end{equation}

\noindent so $[L_{-j}, \, S_r] \neq 0,$ since by Lemma \ref{Lem:2.31} we have

\begin{equation} 0 \neq [[L_{-1}, \, L_1], \, S_r] =  [[L_{-j}, \, S_j], \, S_r] =  [[L_{-j}, \, S_r], \, S_j]\nonumber\end{equation}

\noindent Consequently,

\begin{equation}\label{eq:kleqqj} k \leqq j.\end{equation}

	On the other hand, by definition of $k$, $[L_{-k}, \, S_r] \neq 0$, and $[L_{-i}, \, S_r] = 0$ for $k - 1 \geqq i \geqq 2.$  It follows that if $[L_{-k}, \, S_{k-1}] = 0$ and $[L_{-k}, \, S_k] = 0$, then

\begin{equation} U \eqdef \sum_{i\geqq 0}(\ad L_{-1})^i[L_{-k}, \, S_r]\nonumber\end{equation}

\noindent would be a proper ideal of $S$, to contradict the simplicity of $S$.  Thus, it must be that either $[L_{-k}, \, S_k]  \neq 0$ or $[L_{-k}, \, S_{k-1}] \neq 0$.

	Now, by \eqref{eq:6.1.1} and \eqref{eq:6.3.1}, $j$ is the smallest integer greater than or equal to two such that $[L_{-j}, \, S_j] \neq 0.$  Consequently, if $[L_{-k}, \, S_k] \neq 0,$ then $j \leqq k,$ so by \eqref{eq:kleqqj}, $j=k$.   If, on the other hand, $[L_{-k}, \, S_k] = 0$, then, as we noted above, we must have $[L_{-k}, \, S_{k-1}] \neq 0.$  But by \eqref{eq:6.1.2} and \eqref{eq:6.3.1}, $j+1$ is the smallest $i \geqq 2$ such that $[L_{-i}, \, S_{i-1}] \neq 0,$ so we must have $j+1\leqq k,$ so $j \leqq k-1$, to contradict \eqref{eq:kleqqj}. Thus, 

\begin{equation} j = k\nonumber\end{equation}

\noindent as required. \qed

Define $H$ to be the Lie algebra generated by $L_{-1} $ $ \oplus $ $ L_0 $ $ \oplus $ $ S_1 $ $ \oplus $ $ \dots $ $ \oplus $ $S_{j-1}.$  Since $[L_{-2}, \,  S_1 \oplus$ $\dots \oplus S_{j-1}] = 0,$ it follows from \eqref{eq:6.1} that if $M(H)$ is as in the statement of Weisfeiler's Theorem (Theorem \ref{Thm:1.3}), then $M(H)$ $= L_{-q} \oplus \dots \oplus L_{-2},$ and $H/M(H)$ is a depth-one graded Lie algebra which inherits transitivity (C) and irreducibility (B) from $L$.  From Proposition \ref{Pro:1.4}, we conclude that 
$H/M(H)$ is either between $H(2: \, (1, \, 1),\,
\omega)$ and $CH(2: \, (1, \, 1),\, \omega),$ or is equal to
$L(\epsilon)$ or $M.$  In each of these
cases, $S_1$ is an irreducible abelian $L_0-$module. Thus, we can
from now on assume (See  Corollary \ref{Cor:1.5}.) that assumption (iv) of Section 2 holds and

\begin{equation}\label{eq:6.01} [L_1, \, L_1] = 0 \end{equation}

\noindent and, as a Lie algebra, $[L_{-1}, \, L_1]$ lies between ${\mathfrak{sl}}(2)$ and ${\mathfrak{gl}}(2)$:

\begin{equation}\label{eq:sorg}{\mathfrak{sl}}(2) \subseteq [L_{-1}, \, L_1] \subseteq {\mathfrak{gl}}(2)\end{equation} 	

\begin{Lem}  \label{Lem:2.315} If $[L_{-2}, \, L_1] = 0$ and $[L_{-2}, \, L_2] \neq 0$ and (vii) and \eqref{eq:sorg} hold, then we may assume in what follows that $[L_{-2}, \, L_2] = [L_{-1}, \, L_1].$
\end{Lem}

{\bf{Proof}}  Let $i$ be any integer not equivalent to zero modulo three, and let $z$ be any element of the center of $L_0.$  By transitivity (C), $[L_{-1}, \, z] \neq 0,$ and by irreducibility (B) and Schur's Lemma, $\ad_{L_{-1}}z$ acts as multiplication by a non-zero scalar $a$, so that by (D),  $\ad_{L_{-i}}z$ acts as multiplication by the non-zero scalar $ai$. Consequently, in view of  \eqref{eq:sorg},  Lemma \ref{Lem:1.2}, and transitivity (C) (to deal with $\ad \, z$ when $i$ is positive),

\begin{equation}\label{eq:nonullann} \Ann_{[L_{-1}, \, L_1]}M_{i} = 0 \end{equation}

\noindent for any non-zero $L_0-$submodule $M_i$ of $L_i,  i \not \equiv 0 \,\hbox{(mod 3)}.$    	Since by (D), $[L_{-2}, \, L_2]$ is an ideal of $[L_{-1}, \, L_1]$, we have that $[L_{-2},$ $[L_{-2}, \, L_2]]$ $\neq$ $0$. Similarly,  (See Corollary \ref{Cor:1.41}.)  $[L_2, \, [L_{-2}, L_2]]$ $\neq 0.$  It follows that 

\begin{equation}\label{1grdspneq0} B(L_{-2})_1 \neq 0 \hbox{ and, if $L_2$ is $L_0-$irreducible, } B(L_2)_1 \neq 0.\end{equation}

 	Let $M(B(L_{-2}))$ be as in Weisfeiler's Theorem (Theorem \ref{Thm:1.3}).  We focus on 

\begin{equation} X \eqdef B(L_{-2})/M(B(L_{-2})),\nonumber\end{equation}  

\noindent whose depth is no greater than half of that of $L$.   If the depth and height of $X$ are both greater than one, then we can apply the Main Theorem to conclude that the representation of $L_0'$ on the minus-one component (namely, $L_{-2}$) of $X$ is restricted, so that (See Lemma \ref{Lem:1.2}.) the representation of $L_0'$ on $L_{-1}$ is also restricted, to contradict assumption (vii).  

	If the depth of  $X$ is one, then by Proposition \ref{Pro:1.4},  $X$  is isomorphic either to a Hamiltonian Lie algebra, or to $M,$ or to $L(\epsilon)$ for some $\epsilon.$  Then,  by \eqref{eq:sorg} and (D),

\begin{equation} {\mathfrak{sl}}(2) \subseteq [L_{-2}, \, L_2] \subseteq [L_{-1}, \, L_1] \subseteq {\mathfrak{gl}}(2).\nonumber\end{equation}

\noindent Consequently, we are done unless $[L_{-2}, \, L_2] \cong {\mathfrak{sl}}(2),$ and  $[L_{-1}, \, L_1] \cong {\mathfrak{gl}}(2).$  In that case, we have by Corollary \ref{Cor:1.5} that $X$ is Hamiltonian.  Let $e,$ $f,$ and $h$ be the usual basis of $L_0'.$ Then, according to (1.4) of \cite{BKK}, $ (\ad_{L_{-2}}f)^3$ acts as the identity on the minus-one component of $X,$ namely, $L_{-2}.$   On the other hand, if $\tilde{L}$ is as in Lemma \ref{Lem:2.25}, then $L_{-2} \subseteq M(\tilde{L}),$ and, again by Corollary \ref{Cor:1.5},  either $\tilde{L}/M(\tilde{L}) \cong L(\epsilon),$ or $\tilde{L}/M(\tilde{L}) \cong M.$  If $\tilde{L}/M(\tilde{L}) \cong L(\epsilon)$, then $ (\ad_{L_{-1}}f)^3$ acts as  zero on $L_{-1}.$  On the other hand,  if $\tilde{L}/M(\tilde{L}) \cong M$, then $ (\ad_{L_{-1}}f)^3$ acts as  the identity on $L_{-1}.$  However, by (D) and the fact that  $(\ad_{L_{-1}}f)^3$ is a derivation, we should have that $(\ad_{L_{-2}}f)^3$=$2(\ad_{L_{-1}}f)^3.$  These contradictions enable us to conclude that the lemma is true when the depth of $X$ is one.

	Suppose, finally, that the height of $X$ is one.  Since, as we saw above, $[L_2, \, [L_{-2}, L_2]]$ $\neq 0,$  it follows that $X$ is a non-degenerate Lie algebra. We may therefore apply Corollary \ref{Cor:1.41} to $X$ to conclude that for any irreducible $L_0-$submodule $M_2$ of $L_2$ we have $[M_2, \, [M_2, \, L_{-2}]] \neq 0.$   We apply Corollary \ref{Cor:1.5} to $B(M_{2})/M(B(L_{2})).$ Arguing as in the above paragraph, we obtain a contradiction also in this case.  $~~\square$

\begin{Lem}  \label{Lem:6.021} Suppose $j > 2.$  Then $j$ cannot equal $r,$ so $L_{j+1} \neq 0.$  Similarly, $L_{j+2} \neq 0,$   \end{Lem}

{\bf Proof} Suppose $j = r.$  Then by \eqref{eq:6.3.1} $0 \neq [L_{-j}, \, S_j] = [L_{-r}, \, S_r].$  Then we would have $[[L_{-r}, \, S_r], S_1] = [[L_{-r}, \, S_1], S_r] \subseteq [L_{-(r-1)}, \, S_r] = 0$ by the definition of $k.$  Since, again,  $j = r,$ however, we would have $[L_{-j+1}, \, L_j] = 0,$ to contradict \eqref{eq:6.3.1}.  Thus, $j \neq r.$ 
	
	Now suppose that $j = r - 1.$   Since $j > 2,$ it follows that $r > 3.$  By Lemma \ref{Lem:6.011}, $[L_{-r}, \, L_r] \neq 0.$  Then by  ((D) and) \eqref{eq:nonullann}, $[L_1, \, [L_{-r}, \, L_r]] \neq 0.$  By (iv) (See the discussion before \eqref{eq:6.01}.), we have $L_1 = [L_{-r+1}, \, L_r],$ so by definition of $k$, $[L_1, \, L_i] = 0, \, 1 \leqq i \leqq k - 2 = j - 2 = r - 3.$  If $r > 4,$ then we would have, for example, $0 = [L_{-1}, \, [L_1, \, L_2]] = [[L_{-1}, \, L_1], \, L_2]],$ to contradict \eqref{eq:nonullann}. If $r=4,$ then $j=3,$ so by \eqref{eq:6.3.1}, $[L_{-3}, \, L_3] \neq 0,$ and by \eqref{eq:6.1.3}, $[L_{-3}, \, L_2] = 0.$  Thus, we would have $0 = [L_{-3}, \, 0] = [L_{-3}, \, [L_2, \, L_3]] = [L_2, \, [L_{-3}, \, L_3]],$ to again contradict \eqref{eq:nonullann} by (D).  $~~\square$

Now, by Lemma \ref{Lem:2.25}, $j \leq 2$ in the non-Hamiltonian cases.  Thus, suppose that we are in the Hamiltonian case, and suppose, for a contradiction, that $j > 2.$   For this (Hamiltonian, $j >2$) case, we will assume without loss of generality, that $L$ is generated by $L_{<0} \oplus L_0 \oplus S_1 \oplus \cdots \oplus S_{j-1} \oplus S_j.$  By definition of $j$ and the Jacobi Identity, $[S_i, \, [L_{-2}, \, L_{-2}]] = 0, \, 1 \leqq i \leqq j-1.$  Furthermore, since $j$ is assumed to be greater than two, we have $[S_j, \, [L_{-2}, \, L_{-2}]] \subseteq [L_{-2}, \, S_{j-2}] = 0.$  Thus, again by the Jacobi Identity, $[L_{>0}, \, [L_{-2}, \, L_{-2}]] = 0,$ so, if $[L_{-2}, \, L_{-2}]$ were not equal to zero then $\sum_{i \geqq 0} (\ad L_{-1})^i[L_{-2}, \, L_{-2}]$ would be a proper ideal of $S,$ to contradict the simplicity of $S.$  Thus, we may assume  that

\begin{equation}\label{eq:6.m2a} [L_{-2}, \, L_{-2}] = 0\end{equation}

\noindent In addition, $[L_{-2}, \, S_{j-2}] = 0$ also entails that

\begin{equation}\label{eq:m2s}[L_{-2}, \, [L_{-2}, \, L_i]] = 0, \, 0 \leqq i \leqq 2j - 1\end{equation}

\begin{Lem} \label{Lem:m2m2jp20} If $j > 2,$ then

\begin{equation}[L_{-2}, \, [L_{-2}, \, L_{j+2}]] = 0\nonumber\end{equation}

\noindent\end{Lem}

{\bf Proof.} Since we are assuming that   $j > 2,$ it follows that  $2j > j+2,$ so the lemma follows from \eqref{eq:m2s}. \qed

	We adopt the notation of \cite{BKK}.

\begin{Lem}	\label{Lem:m2jp2ne0} If $j > 2,$ and $j \not\equiv 0$ (mod 3), then $ [L_{-2}, \, [L_{2}, \,  L_{j}]] \neq 0.$ \end{Lem}

{\bf{Proof:}}  Suppose not.  Then

\begin{equation} 0 =  [L_{-2}, \, [L_{2}, \,  L_{j}]] =  [L_{2}, \, [L_{-2}, \,  L_{j}]]\nonumber\end{equation}

\noindent since we are assuming that $j > 2.$  But $0 \neq [L_{-2}, \,  L_{j}] \subseteq H_{j-2}$ (See the discussion preceding \eqref{eq:6.01}.), and $H_{j-2}$ is $L_0-$irreducible.  Consequently, we have $[H_2, \, H_{j-2}] = 0.$  But $\{y^{(3)}, \, xy^{(j-1)}\} = y^{(2)}y^{(j-1)} \equiv jy^{(j+1)}\not\equiv 0$ (since $j \not\equiv 0$ (mod 3)).  \qed

\begin{Lem} \label{Lem:ad1cubed} If $M$ is an $L_0-$submodule of $L$ such that $[L_{-2}, \, M] = 0$ and $[L_{-2}, \, [1, \, M]] $ = $ 0,$ then $(\ad 1)^3$ is an $L_0-$homomorphism of $M$ into $(\ad 1)^3M.$\end{Lem}

{\bf Proof.}  Since $\{1, \, xy\} = 0 = \{1, \, x^{(2)}y\},$ we have  $(\ad 1)^3[L_0,$ $M]$ $=$ $[(\ad 1)^3L_0,$ $M]$ +  $[L_0, \, (\ad 1)^3M]$ $=$ $[(\ad 1)^2x^{(2)}, \, M]$ $+$ $[L_0, \, (\ad 1)^3M]$. Now, $[(\ad 1)^2x^{(2)},$ $ M]$ $\subseteq$ $[1, \, [L_{-2}, \, M]] + [[1, \, M], \, L_{-2}] = 0.$\qed

By transitivity (C), ($S_j \subseteq$)  $L_j$ must be contained in the nine-dimensional $L_0-$module $L_{-1}^*\otimes S_{j-1}.$  (Here, and in what follows, we understand ``is contained in'' in this context to mean ``is $L_0-$isomorphic to an $L_0-$submodule of'' and ``is acted on by $\ad L_{-1}$ as  $\ad L_{-1}$ acts on'' and we will often make use of this implication of (C) without comment.)    Note that if the gradation degree of any of the following three-dimensional $L_0-$modules

\begin{equation} X_{residue\, of\, j\,  (mod \, 3), \, gradation\ degree, \, sequential \ number}\nonumber\end{equation}

\noindent is not equivalent to zero modulo three, then the $L_0-$module is $L_0-$irreducible.

To begin our proof that $j = 1,$ we will assume first that

\begin{equation}\label{eq:jg2} j > 2\end{equation}

	We consider separately the cases where $j$ is equivalent to one, two, or zero modulo three.  We will show that in each case, if $j > 2,$ then $L$ can be replaced with a subalgebra (containing $L_{-1}$ and $L_0$) in which the smallest integer $i > 2$ such that $L_{-2}$ has nonzero bracket with the $i^{th}$ gradation space of the replacement is greater than $j.$

	Note that

\begin{equation} \label{eq:up} \{y^{(2)}, \, y^{(j)}\} = (j+1)^2 x^{(2)}y^{(j+1)} \end{equation}

\noindent Since $L$ contains $S_{j-1} \cong H(2: \mathbf{n}, \omega)_{j-1},$ it follows  that whenever $j \not \equiv 2$ (mod 3),  we also have that $S_j$ contains the non-zero $L_0-$submodule (annihilated by $\ad L_{-2}$) $[S_1, S_{j-1}] \cong H(2: \mathbf{n}, \omega)_{j}.$  Furthermore, if we consider \eqref{eq:up} with $j$ replaced by $j+1,$ we see that if $j \equiv0$ (mod 3), then $S_{j+1}$ contains the non-zero $L_0-$submodule (also annihilated by $\ad L_{-2}$) $[S_1, \, [S_1, \, S_{j-1}]] \cong H(2: \mathbf{n}, \omega)_{j+1}.$

	When $j \equiv 1\, (\hbox{mod 3}), $  $L_{-1}^*\otimes S_{j-1}$$(\supseteq L_j)$ contains the two $L_0-$submodules

\begin{eqnarray}
X_{1, \,  j, \, 1} &\eqdef&  <1^*\otimes x^{(2)}y^{(j)} , \, 2\cdot x^*\otimes x^{(2)}y^{(j)} + 1^*\otimes xy^{(j)},\cr&&  (x^{(2)})^*\otimes x^{(2)}y^{(j)} +  x^*\otimes xy^{(j)} + 1^*\otimes y^{(j)}>             \nonumber\end{eqnarray}

\noindent and 

\begin{eqnarray}
X_{1, \,  j, \, 2} &\eqdef&  <x^*\otimes x^{(2)}y^{(j)} , \, 2\cdot (x^{(2)})^*\otimes x^{(2)}y^{(j)} + x^*\otimes xy^{(j)},\cr&&  1^*\otimes x^{(2)}y^{(j)} +  (x^{(2)})^*\otimes xy^{(j)} + x^*\otimes y^{(j)}>             \nonumber\end{eqnarray}

	(Acting here like $\frac{\partial}{\partial x}$, $\ad y$ maps $xy^{(j)}$,  $x^{(2)}y^{(j)}$, 
and  $y^{(j)}$, to  $y^{(j)}$, $xy^{(j)}$, and  zero, respectively.  Also, adding one to the divided power modulo three and multiplying by two,  $\ad y$ maps $1^*$,  $x^*$, and  $(x^{(2)})^*$, to $2\cdot x^*$, $2\cdot (x^{(2)})^*$ and  $2\cdot 1^*$, respectively.)

	 Now, 

\begin{equation} \label{eq:lj}L_j \subseteq L_{-1}^*\otimes S_{j-1} = X_{1, \,  j, \, 1} \oplus X_{1, \,  j, \, 3},\end{equation} 

\noindent where modulo its submodule  $X_{1, \,  j, \, 2},$ the indecomposable $L_0-$module $X_{1, \,  j, \, 3}$ is spanned by $(x^{(2)})^*\otimes x^{(2)}y^{(j)} , \, 2\cdot 1^*\otimes x^{(2)}y^{(j)} + (x^{(2)})^*\otimes xy^{(j)}, \hbox{ and }  1^*\otimes xy^{(j)} +  x^*\otimes x^{(2)}y^{(j)} + (x^{(2)})^*\otimes y^{(j)}.$

 Now $[L_{-2}, \, X_{1, \,  j, \, 2}] = 0,$  and, as an $L_0-$submodule of $L_{-1}^*\otimes S_{j-1},$ $X_{1, \,  j, \, 2}$ is unique in this regard.  It follows from the discussion following \eqref{eq:up} that $S_j \subseteq L_j$ must contain an $L_0-$submodule 

\begin{equation} \label{eq:ham} Q_j \cong X_{1, \,  j, \, 2} \cong H(2: \mathbf{n}, \omega)_j = \langle  y^{(j+1)}, \, x y^{(j+1)}, \, x^{(2)}y^{(j+1)}\rangle.\end{equation}

\noindent Since $S_j = [L_{-1}, \, S_{j+1}]$ by Lemma \ref{Lem:2.3},  it follows from transitivity (C) that
 $S_{j+1}$ must be contained in

\begin{equation} L_{-1}^*\otimes S_j \subseteq L_{-1}^*\otimes L_j\subseteq L_{-1}^*\otimes (X_{1, \,  j, \, 1} \oplus X_{1, \,  j, \, 3}) =  L_{-1}^*\otimes X_{1, \,  j, \, 1} \oplus  L_{-1}^*\otimes  X_{1, \,  j, \, 3}\nonumber\end{equation}

\noindent Since $Q_j \subseteq S_j,$  $S_{j+1}$ must have non-zero intersection with $L_{-1}^*\otimes Q_j \cong  L_{-1}^*\otimes X_{1, \, j, \, 2}.$   (Note that if $X_{1, \, j, \, 2} \cap S_j = 0,$ then $X_{1, \, j, \, 3} \cap S_j = 0;$   similarly, if $(L_{-1}^*\otimes X_{1, \, j, \, 2}) \cap S_{j+1} = 0,$ then $(L_{-1}^*\otimes X_{1, \, j, \, 3}) \cap S_{j+1} = 0$.)

	If  $S_{j+1} \cap (L_{-1}^*\otimes X_{1, \, j, \, 2})$ were to contain an $L_0-$submodule $R_{j+1}$ such that $[L_{-2}, \, R_{j+1}] \neq 0,$ we could replace $L$ with the Lie algebra generated by $L_{\leqq 0} \oplus S_1 \oplus \cdots \oplus S_{j-1} \oplus Q_j \oplus R_{j+1},$ and we would have a Lie subalgebra of $L$ such that the minimal $i > 2,$ such that $L_{-2}$ has non-zero bracket with the $i^{th}$ gradation space of the subalgebra, is $j+1.$ 

	Define

\begin{eqnarray} d &\eqdef& x^* \otimes x^{(2)}y^{(j)}\cr
e &\eqdef& 2\cdot (x^{(2)})^*\otimes x^{(2)}y^{(j)} + x^*\otimes xy^{(j)}\cr
f &\eqdef& 1^*\otimes x^{(2)}y^{(j)} + (x^{(2)})^*\otimes xy^{(j)} + x^*\otimes y^{(j)}\cr\nonumber\end{eqnarray}

\noindent Then $y \cdot f = 2d,$ and $L_{-1}^*\otimes X_{1, \, j, \, 2}$ is the sum of the following $L_0-$submodules:

\begin{equation} X_{1, \, j+1, \, 1} \eqdef \langle (x^{(2)})^*\otimes e + x^*\otimes f + 2\cdot  1^*\otimes d\rangle + \langle 2\cdot x^*\otimes e + 1^*\otimes f, \,  1^*\otimes e\rangle\nonumber\end{equation}

\begin{equation} X_{1, \, j+1, \, 2} \eqdef \langle (x^{(2)})^*\otimes e + x^*\otimes f +  1^*\otimes d\rangle + \langle x^*\otimes e + 2\cdot (x^{(2)})^*\otimes d, \,  x^*\otimes d\rangle\nonumber\end{equation}

\begin{equation} X_{1, \, j+1, \, 3} \eqdef \langle (x^{(2)})^*\otimes f + 2 \cdot x^*\otimes d + 1^*\otimes e, \,  x^*\otimes f + 2\cdot (x^{(2)})^*\otimes e, \,  x^*\otimes e\rangle\nonumber\end{equation}

\noindent Note that 

\begin{equation}\label{eq:z1jp16}[L_{-2}, \,  X_{1, \, j+1, \, 2}] = 0.\end{equation}

	Now, since $(\ad [1, \, x])\cdot (1^* \otimes e ) = [e, \, x] = 2\cdot[x, \, e] = 2xy^{(j)} \neq 0,$ and  $(\ad [1, \, x])\cdot ((x^{(2)})^*\otimes e + x^*\otimes f + 2\cdot  1^*\otimes d) = [1, \, f] + 2\cdot [d, \, x] = 2x^{(2)}y^{(j)}\neq 0 ,$ if $L_{j+1} \cap X_{1, \, j+1, \, 1} \neq 0,$ then $[L_{-2}, \, L_{j+1} \cap X_{1, \, j+1, \, 1}]\neq 0.$ Thus, if we set $R_{j+1} = L_{j+1} \cap X_{1, \, j+1, \, 1}$  in the above argument, we see that we can assume in what follows that

\begin{equation} \label{eq:x1jp15} L_{j+1} \cap X_{1, \, j+1, \, 1} = 0.\end{equation}

	In addition, since $(\ad [x, \, x^{(2)}])\cdot (x^* \otimes e ) = [e, \, x^{(2)}] = 2\cdot [x^{(2)}, \, e] = x^{(2)}y^{(j)} \neq 0,$ and  $(\ad [1, \, x])\cdot (2\cdot (x^{(2)})^*\otimes e + x^*\otimes f) = [1, \, f] = x^{(2)}y^{(j)}\neq 0 ,$ if $L_{j+1} \cap X_{1, \, j+1, \, 3} \neq 0,$ then $[L_{-2}, \, L_{j+1} \cap X_{1, \, j+1, \, 3}]\neq 0.$  Thus, if we set $R_{j+1} = L_{j+1} \cap X_{1, \, j+1, \, 3}$  in the above argument, we see that we can assume in what follows that

\begin{equation} \label{eq:x1jp17} L_{j+1} \cap X_{1, \, j+1, \, 3} = 0.\end{equation}

	Now, $X_{1, \, j+1, \, 1}$ and $X_{1, \, j+1, \, 2}$ are isomorphic as $L_0-$modules, so  suppose that, for scalars $A,$ and $B,$ $L_{j+1} \cap (X_{1, \, j+1, \, 1}  + X_{1, \, j+1, \, 2})$ $\supseteq$ $AX_{1, \, j+1, \, 1} + BX_{1, \, j+1, \, 2},$ If this latter module had non-zero bracket with $L_{-2},$  then we could set $R_{j+1}$ equal to it and argue as in the previous two paragraphs.  If, on the other hand, the bracket with $L_{-2}$ were zero, then, for example, we would have that the bracket of $[1, \, x]$ with

\begin{equation}  A\cdot ((x^{(2)})^*\otimes e + x^*\otimes f + 2\cdot 1^*\otimes d) + B\cdot((x^{(2)})^*\otimes e + x^*\otimes f +  1^*\otimes d)\nonumber\end{equation}

\noindent was zero.  Consequently, we would have

\begin{equation} 0 = A\cdot [1, \, f] + 2A\cdot [d, \, x] + B\cdot [1, \, f] + B\cdot [d, \, x] = 2Ax^{(2)}y^{(j)} \nonumber\end{equation}

\noindent from which we can conclude that $A = 0,$   so that the only linear combination of $X_{1, \, j+1, \, 1}$ and $X_{1, \, j+1, \, 2}$ to have zero bracket with $L_{-2}$ is $X_{1, \, j+1, \, 2}.$  It follows that we may assume (See \eqref{eq:up}ff.) that $L_{j+1}$ has non-zero intersection with $X_{1, \, j+1, \, 2}.$

	Define

\begin{eqnarray} \kappa &\eqdef& x^*\otimes d\cr
\lambda &\eqdef&  x^*\otimes e + 2\cdot (x^{(2)})^*\otimes d\cr
\mu &\eqdef& (x^{(2)})^*\otimes e + x^*\otimes f + 1^*\otimes d\nonumber\end{eqnarray}

\noindent Then $y \cdot \mu = \kappa,$ and $L_{-1}^*\otimes X_{1, \, j+1, \, 2},$ which by transitivity (C) has non-zero intersection with $L_{j+2},$  is the sum of the following (non-isomorphic) $L_0-$modules:

\begin{equation} X_{1, \, j+2, \, 1} \eqdef <1^*\otimes \mu + (x^{(2)})^*\otimes \kappa + x^*\otimes \lambda> + <2\cdot (x^{(2)})^*\otimes \mu + x^*\otimes \kappa, \, x^*\otimes \mu>\nonumber\end{equation}

\begin{equation} X_{1, \, j+2, \, 2} \eqdef <1^*\otimes \kappa +  (x^{(2)})^*\otimes \lambda + x^*\otimes \mu, \, 2\cdot (x^{(2)})^*\otimes \kappa + x^*\otimes \lambda, \, x^*\otimes \kappa> \nonumber\end{equation}

\begin{equation}X_{1, \, j+2, \, 3} \eqdef < (x^{(2)})^*\otimes \mu + x^*\otimes \kappa + 1^*\otimes \lambda, \, 2\cdot x^*\otimes \mu + 1^* \otimes \kappa> + <1^*\otimes \mu>\nonumber\end{equation}

	Now, since $(\ad [1, \, x])\cdot (x^* \otimes \mu ) = [1, \, \mu] =  d,$ and  $(\ad [1, \, x])\cdot ( (x^{(2)})^*\otimes \kappa +   x^*\otimes \lambda + 1^*\otimes \mu) = [\mu, \, x] + [1, \, \lambda] = 2f\neq 0 ,$ if $L_{j+2} \cap X_{1, \, j+2, \, 1} \neq 0,$ then $[L_{-2}, \, L_{j+2} \cap X_{1, \, j+2, \, 1}]\neq 0,$ and we can replace $L$ with the Lie algebra generated by
$L_{\leqq 0} \oplus L_1 \oplus \cdots \oplus L_{j-1} \oplus (L_j \cap X_{1, \, j, \, 2}) \oplus  (L_{j+1} \cap X_{1, \, j+1, \, 2}) \oplus (L_{j+2} \cap X_{1, \, j+2, \, 1}),$ and we will have a Lie subalgebra of $L$ such that the minimal $i > 2,$ such that $L_{-2}$ has non-zero bracket with the $i^{th}$ gradation space of the subalgebra, is $j+2.$    Thus, we may assume in what follows that

\begin{equation} \label{eq:x114} L_{j+2} \cap X_{1, \, j+2, \, 1} = 0.\end{equation}

	Also, since $(\ad [1, \, x])\cdot (1^* \otimes \mu ) = [\mu, \, x] =  2f,$ and  $(\ad [1, \, x])\cdot (2\cdot x^*\otimes \mu + 1^*\otimes \kappa) = 2\cdot [1, \, \mu] + [\kappa, \, x] = 2d + 2d \not\equiv 0 ,$ if $L_{j+2} \cap X_{1, \, j+2, \, 3} \neq 0,$ then $[L_{-2}, \, L_{j+2} \cap X_{1, \, j+2, \, 3}]\neq 0,$ and we can replace $L$ with the Lie algebra generated by
$L_{\leqq 0} \oplus L_1 \oplus \cdots \oplus L_{j-1} \oplus (L_j \cap X_{1, \, j, \, 2}) \oplus   (L_{j+1} \cap X_{1, \, j+1, \, 2}) \oplus (L_{j+2} \cap X_{1, \, j+2, \, 3}),$ and we will have a Lie subalgebra of $L$ such that the minimal $i > 2,$ such that $L_{-2}$ has non-zero bracket with the $i^{th}$ gradation space of the subalgebra, is $j+2.$    Thus, we may assume in what follows that

\begin{equation} \label{eq:x116} L_{j+2} \cap X_{1, \, j+2, \, 3} = 0.\end{equation}

	Note that

\begin{equation}\label{eq:z115} [L_{-2}, \, X_{1, \, j+2, \, 2}] = 0,\end{equation}

\noindent from which we can conclude that if $L_{j+2}$ has non-zero intersection with $X_{1, \, j+2, \, 2}$, then by \eqref{eq:z1jp16} and Lemma \ref{Lem:ad1cubed}, we have that $L_{j+2} \cap X_{1, \, j+2, \, 2} \cong S_{j-1}.$  We could then repeat the above argument (with $S_{j-1}$ replaced by $L_{j+2} \cap X_{1, \, j+2, \, 2}$) to either obtain a Lie algebra in which the minimal gradation degree $i > 2$ such that $[L_{-2}, \, L_i] \neq 0$ is each time even greater, or,  eventually, to find that the highest gradation space $S_s$ of $S$ contains an $L_0-$submodule $Q_s$ which has zero bracket with $L_{-2}.$  However, by Lemma \ref{Lem:2.3}, $S_s$ is irreducible as an $L_0-$module, so we would have that $[L_{-2}, \, S_s] = 0.$  Since $S$ is generated by $S_s$ and $L_{-1},$ it would follow from the construction of $Q_s$ $ (= S_s)$ (See also Lemma \ref{Lem:2.3}.) that $[L_{-2}, \, S_{>0}] = 0,$ to contradict, for example, \eqref{eq:6.00001}.

	In view of  \eqref{eq:x114}, \eqref{eq:x116},  and  \eqref{eq:z115}ff,   if $j \equiv 1$ (mod 3), we have shown (since, from above, $L_{-1}^*\otimes X_{1, \, j+1, \, 2}$ must have non-zero intersection with $L_{j+2}$)  that we can assume that $L$ contains a Lie subalgebra such that  if $i > 2$ is minimal such that $L_{-2}$ has non-zero bracket with the $i^{th}$ gradation space of that subalgebra, then $i$  is greater than $j.$

	If $j \equiv 2$ (mod 3), then $\{y, \, y^{(j)}\} = 2x^{(2)}y^{(j)},$ and $L_j \subseteq (L_{-1})^* \otimes S_{j-1},$ which equals the sum of the $L_0-$submodules

\begin{eqnarray}X_{2, \,  j, \, 1} &\eqdef&  <1^*\otimes xy^{(j)} , \, 2\cdot x^*\otimes xy^{(j)} + 1^*\otimes y^{(j)},\cr&&  (x^{(2)})^*\otimes xy^{(j)} +  x^*\otimes y^{(j)} + 2\cdot1^*\otimes x^{(2)}y^{(j)}>             \nonumber\end{eqnarray}

\noindent and 

\begin{eqnarray}
X_{2, \,  j, \, 2} &\eqdef&  <x^*\otimes x^{(2)}y^{(j)} , \, 2\cdot (x^{(2)})^*\otimes x^{(2)}y^{(j)} + x^*\otimes xy^{(j)},\cr&&  1^*\otimes x^{(2)}y^{(j)} +  (x^{(2)})^*\otimes xy^{(j)} + x^*\otimes y^{(j)}>             \nonumber\end{eqnarray}

\noindent and

\begin{eqnarray}
X_{2, \,  j, \, 3} &\eqdef&  <x^*\otimes xy^{(j)} , \, 2\cdot  (x^{(2)})^*\otimes xy^{(j)} + x^*\otimes y^{(j)},\cr&&  (x^{(2)})^*\otimes y^{(j)} +  2\cdot x^*\otimes  x^{(2)}y^{(j)} + 1^*\otimes xy^{(j)}> \nonumber\end{eqnarray}

	Since $[L_{-2}, \, X_{2, \,  j, \, 2}] = 0,$ it follows from the definition of $j$ \eqref{eq:6.00001} that $S_j \not \subset  X_{2, \, j, \, 2},$ so either $X_{2, \,  j, \, 1} \cap S_j \neq 0$ or $X_{2, \,  j, \, 3} \cap S_j \neq 0.$

	Let us first suppose that $S_j \cap X_{2, \,  j, \, 1} \neq 0$.  Define
 
\begin{eqnarray} a &\eqdef& 1^*\otimes xy^{(j)}\cr
b &\eqdef& 2\cdot x^*\otimes xy^{(j)} + 1^* \otimes y^{(j)} \hbox{ and }\cr
c &\eqdef& (x^{(2)})^*\otimes xy^{(j)} + x^*\otimes y^{(j)} + 2 \cdot 1^*\otimes x^{(2)}y^{(j)}\nonumber\end{eqnarray}

Then $y \cdot c = a,$ and by transitivity (C), $L_{j+1}$ has non-zero intersection with  $ L_{-1}^* \otimes X_{2, \,  j, \, 1},$ which  is the sum of the following (non-isomorphic) $L_0-$submodules:

\begin{equation} X_{2, \,  j+1, \, 1} \eqdef \langle (x^{(2)})^*\otimes a + x^* \otimes b + 1^* \otimes c\rangle + \langle 2\cdot x^*\otimes a + 1^* \otimes b,  \, 1^* \otimes a\rangle\nonumber\end{equation}

\begin{equation} X_{2, \,  j+1, \, 2} \eqdef  \langle  (x^{(2)})^*\otimes c + x^* \otimes a + 1^* \otimes b, \, 2\cdot x^*\otimes c + 1^* \otimes a \rangle +  \langle1^* \otimes c\rangle \nonumber\end{equation}

\begin{equation} X_{2, \,  j+1, \, 3} \eqdef  \langle  (x^{(2)})^*\otimes b + x^* \otimes c + 1^* \otimes a, \, 2\cdot (x^{(2)})^*\otimes a + x^* \otimes b, \,  x^* \otimes a\rangle \nonumber\end{equation}

	 We  first suppose that $S_{j+1}$ has non-zero intersection with $X_{2, \,  j+1, \, 1},$  and we set

\begin{eqnarray}\alpha &\eqdef& 1^*\otimes a \cr
\beta &\eqdef& 2\cdot x^*\otimes a + 1^*\otimes b\cr
\gamma &\eqdef& (x^{(2)})^*\otimes a + x^*\otimes b + 1^*\otimes c\nonumber\end{eqnarray}

\noindent Then $y \cdot \gamma = \alpha,$ and by transitivity (C), $S_{j+2}$ would have non-zero intersection with  $L_{-1}^*\otimes X_{2, \,  j+1, \, 1},$ which is the sum of the following $L_0-$submodules:

\begin{equation} X_{2, \,  j+2, \, 1} \eqdef \langle (x^{(2)})^*\otimes \gamma, \, 2\cdot x^*\otimes \gamma\rangle, + \langle 1^*\otimes \gamma \rangle \nonumber\end{equation}

\begin{equation} X_{2, \,  j+2, \, 2} \eqdef \langle (x^{(2)})^*\otimes \alpha + x^*\otimes \beta +  1^*\otimes \gamma\rangle + \langle 2\cdot x^*\otimes \alpha +  1^*\otimes \beta, \, 1^*\otimes \alpha \rangle \nonumber\end{equation}

\begin{equation} X_{2, \,  j+2, \, 3} \eqdef \langle (x^{(2)})^*\otimes \beta + x^*\otimes \gamma +  1^*\otimes \alpha, \, x^*\otimes \beta +  2\cdot (x^{(2)})^*\otimes \alpha, \, x^*\otimes \alpha \rangle \nonumber\end{equation}

	Focusing on $X_{2, \,  j+2, \, 1},$ we have $[1, \, x]\cdot (1^*\otimes \gamma) = [\gamma, \, x] = 2\cdot [x, \, \gamma] = 2b,$ and $(\ad [1, \, x]) \cdot b =  (\ad [1, \, x])\cdot (2\cdot x^*\otimes xy^{(j)} + 1^*\otimes y^{(j)}) = 2\cdot [1, \, xy^{(j)}] + [y^{(j)}, \, x] = y^{(j-1)} \neq 0,$
so by Lemma \ref{Lem:m2m2jp20},  

\begin{equation}\label{eq:x29}X_{2, \,  j+2, \, 1} \cap L_{j+2} = 0.\end{equation}

	Focusing on $X_{2, \,  j+2, \, 3},$ we have $(\ad [1, \, x])\cdot (x^*\otimes \alpha) = [1, \, \alpha] = a,$ and $(\ad [1, \, x]) \cdot a =  (\ad [1, \, x])\cdot (1^*\otimes xy^{(j)}) = [xy^{(j)}, \, x] = xy^{(j-1)} \neq 0,$
so by Lemma \ref{Lem:m2m2jp20} again,  

\begin{equation}\label{eq:x211}X_{2, \,  j+2, \, 3} \cap L_{j+2} = 0.\end{equation}

	It follows that if $ X_{2, \,  j+1, \, 1} \subset S_{j+1},$ then $ X_{2, \,  j+2, \, 2} \subset L_{j+2}.$  	Note that 

\begin{equation}\label{eq:z210}[L_{-2}, \, X_{2, \,  j+2, \, 2}] = 0.\end{equation}

\noindent Since also $[L_{-2}, \, X_{2, \,  j+1, \, 1}] = 0,$ it follows that if $L_{j+2} \cap X_{2, \,  j+2, \, 2} \neq 0,$ then we have by Lemma \ref{Lem:ad1cubed} that $(\ad 1)^3$ is an isomorphism from $L_{j+2} \cap X_{2, \,  j+2, \, 2}$ to $S_{j-1}:$

\begin{eqnarray} &(\ad 1)^3\cdot (1^*\otimes \alpha) = (\ad 1)^2\alpha = (\ad 1)a = xy^{(j)}\cr
&(\ad 1)^3\cdot (2\cdot x^*\otimes \alpha + 1^*\otimes \beta) = (\ad 1)^2\beta = (\ad 1)b = y^{(j)}\cr
&(\ad 1)^3\cdot ((x^{(2)})^*\otimes \alpha + x^*\otimes \beta + 1^*\otimes \gamma) = (\ad 1)^2\gamma = (\ad 1)c = 2x^{(2)}y^{(j)}
\nonumber\end{eqnarray}

\noindent  Since $X_{2, \,  j+2, \, 2} \cong S_{j-1}$ as an $L_0-$module, and since the construction of the \lq\lq$X$\rq\rq modules depends only on the $L_0-$module properties of $S_{j-1}$ and $L_{-1}^*,$ and $L_{-1}$ (and its bracket with itself, $L_{-2}$) only interact with the first factor of the tensors, we can repeat the process with $X_{2, \,  j+2, \, 2}$  in place of $S_{j-1}$ and continue to repeat it until we arrive at the highest gradation space. Since $(\ad [1, \, x])\cdot (1^*\otimes xy^{(j)}) = xy^{(j-1)} \neq 0,$ it follows that $[L_{-2}, \, X_{2,j,1}] \neq 0.$  Consequently, when we arrive at the highest gradation space, $L_{-2}$ will have non-zero bracket with either $L_{r-2}$, or $L_{r-1}$, or $L_{r}$.  Then by Lemma \ref{Lem:6.02}, we would have

\begin{equation}\label{eq:jleqq2a}j = k \leqq 2.\end{equation}

\noindent Since \eqref{eq:jleqq2a} is what we are presently trying to establish, we may assume in what follows that $L_{j+2} \cap X_{2, \,  j+2, \, 2} = 0,$ which equation, together with \eqref{eq:x29} and \eqref{eq:x211}, shows that $L_{-1}^*\otimes X_{2, \,  j+1, \, 1} \cap L_{j+2} = 0,$ so that we may assume in what follows that $L_{j+1}$ has zero intersection with $X_{2, \,  j+1, \, 1}.$

	Next, suppose that  $X_{2, \,  j+1, \, 2}  \cap L_{j+1} \neq 0.$   Then $ \tilde{X}_{2, \,  j+1, \, 2}  \cap L_{j+1} \neq 0,$ where

\begin{equation} \tilde{X}_{2, \,  j+1, \, 2} \eqdef \langle  (x^{(2)})^*\otimes c + x^* \otimes a + 1^* \otimes b, \, 2\cdot x^*\otimes c + 1^* \otimes a \rangle,\nonumber\end{equation}

\noindent since  $\tilde{X}_{2, \,  j+1, \, 2}$ is the irreducible $L_0-$submodule of the indecomposable $L_0-$module $ X_{2, \,  j+1, \, 2}.$  Set

\begin{eqnarray} \delta &\eqdef&  (x^{(2)})^*\otimes c + x^* \otimes a + 1^* \otimes b \hbox{ and }\cr
\epsilon &\eqdef& 2\cdot x^*\otimes c + 1^* \otimes a\nonumber\end{eqnarray}

\noindent Then  $y \cdot \delta = 0,$ and $ L_{-1}^* \otimes \tilde{X}_{2, \,  j+1, \, 2}$ is the sum of the following irreducible (non-isomorphic) $L_0-$submodules

\begin{equation}X_{2, \,  j+2, \, 4}\eqdef  \langle (x^{(2)})^*\otimes\epsilon+ x^*\otimes \delta\rangle +  \langle 2\cdot x^*\otimes\epsilon+ 1^*\otimes \delta, 1^* \otimes \epsilon\rangle\nonumber\end{equation}

\noindent and

\begin{equation} X_{2, \,  j+2, \, 5} \eqdef \langle (x^{(2)})^*\otimes \delta +  1^*\otimes \epsilon,  2\cdot (x^{(2)})^*\otimes\epsilon+ x^*\otimes \delta, x^* \otimes \epsilon\rangle\nonumber\end{equation}

	If, in addition, we set

\begin{equation}\zeta \eqdef 1^*\otimes c\nonumber\end{equation}

\noindent then $y \cdot \zeta = \epsilon,$ and

\begin{equation}L_{-1}^*\otimes X_{2, \, j+1, \, 2} =  L_{-1}^*\otimes \tilde{X}_{2, \,  j+1, \, 2} + X_{2, \, j+2, \, 6} =   X_{2, \,  j+2, \, 4} + X_{2, \,  j+2, \, 5} +  X_{2, \, j+2, \, 6}\nonumber\end{equation}

\noindent where
\begin{equation} X_{2, \, j+2, \, 6} \eqdef <1^* \otimes \zeta> + <2\cdot x^* \otimes \zeta + 1^* \otimes \epsilon, \, (x^{(2)})^*\otimes \zeta + x^* \otimes\epsilon+ 1^*\otimes \delta>\nonumber\end{equation}

\noindent Note that 

\begin{equation}\label{eq:z218} [L_{-2}, \, X_{2, \, j+2, \, 6}] = 0.\end{equation}

	Consider $X_{2, \,  j+2, \, 4}.$  We have

\begin{equation}(\ad [1, \, x])\cdot(1^*\otimes \epsilon) = [\epsilon, \, x] = 2\cdot [x, \, \epsilon] = c\nonumber\end{equation}

\noindent and

\begin{equation}(\ad [1, \, x])\cdot c = 2\cdot [x^{(2)}y^{(j)}, \, x] + [1, \, y^{(j)}] = 2x^{(2)}y^{(j-1)} - x^{(2)}y^{(j-1)} = x^{(2)}y^{(j-1)} \neq 0\nonumber\end{equation}

\noindent so by Lemma \ref{Lem:m2m2jp20}, 

\begin{equation}\label{eq:x21}X_{2, \,  j+2, \, 4} \cap L_{j+2} = 0.\end{equation}

Next consider $X_{2, \,  j+2, \, 5}.$  Here we have

\begin{equation} (\ad [1, \, x])\cdot (x^*\otimes \epsilon) = [1, \, \epsilon] = a\nonumber\end{equation}

\noindent and (from above) $(\ad [1, \, x])\cdot a \neq 0,$ so by Lemma \ref{Lem:m2m2jp20} again,  

\begin{equation}\label{eq:x22}X_{2, \,  j+2, \, 5} \cap  L_{j+2} = 0.\end{equation}

\noindent  It follows that 

\begin{equation}\tilde{X}_{2, \,  j+1, \, 2}  \cap L_{j+1} = 0.\end{equation} 

\noindent This, from above, implies that

\begin{equation}\label{eq:x212}X_{2, \,  j+1, \, 2}  \cap L_{j+1} = 0.\end{equation}

	We next suppose that  $X_{2, \,  j+1, \, 3} \cap L_{j+1} \neq 0$. We set

\begin{eqnarray}\eta &\eqdef&  x^*\otimes a\cr
\theta &\eqdef& 2 \cdot  (x^{(2)})^*\otimes a + x^* \otimes b\cr
\iota &\eqdef& (x^{(2)})^*\otimes b + x^*\otimes c +  1^*\otimes a\nonumber\end{eqnarray}

\noindent Then $y \cdot \iota = 0,$ and $L_{-1}^*\otimes X_{2, \,  j+1, \, 3}$ contains the following two (non-isomorphic)  $L_0-$submodules:

\begin{equation}X_{2, \,  j+2, \, 7}\eqdef  \langle (x^{(2)})^*\otimes \eta + x^*\otimes \theta + 1^* \otimes \iota, \,                              2\cdot x^*\otimes \eta + 1^*\otimes \theta, 1^* \otimes \eta\rangle\nonumber\end{equation}

\noindent and

\begin{equation} X_{2, \,  j+2, \, 8} \eqdef \langle (x^{(2)})^*\otimes \theta + x^* \otimes \iota + 1^*\otimes \eta,  2\cdot (x^{(2)})^*\otimes \eta + x^*\otimes \theta, x^* \otimes \eta\rangle\nonumber\end{equation}

\noindent Furthermore, $L_{-1}^*\otimes X_{2, \,  j+1, \, 3} = X_{2, \,  j+2, \, 7}  + X_{2, \,  j+2, \, 9},$ where, modulo its irreducible submodule $ X_{2, \,  j+2, \, 8},$ $X_{2, \,  j+2, \, 9}$ is spanned by $(x^{(2)})^*\otimes \eta,$ $(x^{(2)})^*\otimes \theta + 2\cdot1^*\otimes \eta ,$ and $(x^{(2)})^*\otimes \iota + x^*\otimes \eta +  1^*\otimes \theta.$  Now,

\begin{eqnarray} (\ad [x, \, x^{(2)}])\cdot ((x^{(2)})^*\otimes \eta) &=& [x, \, \eta] = a\cr
(\ad [x^{(2)}, \, x])\cdot ((x^{(2)})^*\otimes \theta + 2\cdot1^*\otimes \eta) &=& [\theta, \, x] = 2b\cr
(\ad [1, \, x])\cdot ((x^{(2)})^*\otimes \iota + x^*\otimes \eta +  1^*\otimes \theta) &=& [1, \, \eta] +  [\theta, \, x] = 2b\nonumber\end{eqnarray}

\noindent so, since from above we know that $(\ad [1, \, x])\cdot a$ and $(\ad [1, \, x])\cdot b$ are both non-zero, it follows from Lemma \ref{Lem:m2m2jp20} that 

\begin{equation}\label{eq:x217}L_{j+2} \cap X_{2, \,  j+2, \, 9} \subseteq X_{2, \,  j+2, \, 8},\end{equation}

\noindent Note that 

\begin{equation}\label{eq:z216} [L_{-2}, \, X_{2, \, j+2, \, 8}] = 0.\end{equation}

Focusing on $X_{2, \,  j+2, \, 7}$  we have

\begin{equation} (\ad [1, \, x])\cdot (1^*\otimes \eta) = [\eta, \, x] = 2a\nonumber\end{equation}

\noindent  so, since, again, from above we know that $(\ad [1, \, x])\cdot a$ is not zero, it follows from Lemma \ref{Lem:m2m2jp20} that 

\begin{equation}\label{eq:x215}X_{2, \,  j+2, \, 7} \cap  L_{j+2} = 0.\end{equation}

\noindent It follows from \eqref{eq:x217}, \eqref{eq:z216}, and \eqref{eq:x215} that $L_{-1}^*\otimes X_{2, \,  j+1, \, 3}$ has zero bracket with  $L_{-2}.$

 	Now suppose that $X_{2, \,  j, \, 3} \cap L_j \neq 0.$  Set

\begin{eqnarray} d &\eqdef& x^*\otimes xy^{(j)}\cr
e &\eqdef& 2 \cdot (x^{(2)})^*\otimes xy^{(j)} + x^* \otimes y^{(j)} \hbox{ and }\cr
f &\eqdef& (x^{(2)})^*\otimes y^{(j)} + 2 \cdot x^*\otimes x^{(2)}y^{(j)} +  1^*\otimes xy^{(j)}\nonumber\end{eqnarray}

Then $y \cdot f = d,$ and $ L_{-1}^* \otimes X_{2, \,  j, \, 3}$  contains the following submodules:

\begin{equation}X_{2, \,  j+1, \, 4} \eqdef \langle (x^{(2)})^*\otimes e + x^* \otimes f + 1^* \otimes d,  \, 2\cdot (x^{(2)})^*\otimes d + x^* \otimes  e \rangle  + \langle  x^* \otimes d \rangle\nonumber\end{equation}

\begin{equation} X_{2, \,  j+1, \, 5} \eqdef  \langle  (x^{(2)})^*\otimes d + x^* \otimes e + 1^* \otimes f, \, 2\cdot x^*\otimes d + 1^* \otimes e, \, 1^* \otimes d\rangle \nonumber\end{equation}

\noindent Furthermore, $ L_{-1}^* \otimes X_{2, \,  j, \, 3} = X_{2, \,  j+1, \, 6} + X_{2, \,  j+1, \, 5},$ where, modulo $X_{2, \,  j+1, \, 4},$ $X_{2, \,  j+1, \, 6}$ is spanned by $\langle 1^*\otimes e + x^*\otimes d + (x^{(2)})^*\otimes f\rangle$ and $\langle 2 \cdot (x^{(2)})^*\otimes e + x^*\otimes f, \, x^*\otimes e\rangle.$

	Focusing first on $X_{2, \,  j+1, \, 4},$ we set

\begin{eqnarray} \kappa &\eqdef&  x^* \otimes d\cr
\lambda &\eqdef& 2\cdot (x^{(2)})^*\otimes d + x^* \otimes  e \cr
\mu&\eqdef&  (x^{(2)})^*\otimes e + x^* \otimes f + 1^* \otimes d\nonumber\end{eqnarray}

\noindent Then $y \cdot \mu = 0,$ and $L_{-1}^*\otimes X_{2, \,  j+1, \, 4}$ is the sum of the following three $L_0-$submodules:

\begin{equation} X_{2, \,  j+2, \, 10} \eqdef \langle (x^{(2)})^* \otimes \lambda+ x^*\otimes \mu\rangle+\langle 2\cdot x^*\otimes \lambda + 1^*\otimes \mu, \, 1^*\otimes \lambda\rangle\nonumber\end{equation}

\noindent and

\begin{equation}X_{2, \,  j+2, \, 11} \eqdef \langle (x^{(2)})^*\otimes \kappa + x^*\otimes \lambda + 1^*\otimes \mu, \, 2\cdot x^*\otimes \kappa + 1^*\otimes \lambda\rangle+\langle1^*\otimes \kappa\rangle\nonumber\end{equation}

\noindent and

\begin{equation}X_{2, \,  j+2, \, 12} \eqdef \langle 1^*\otimes \lambda + (x^{(2)})^*\otimes \mu, \, 2\cdot (x^{(2)})^*\otimes \lambda +  x^*\otimes \mu, \, x^*\otimes \lambda\rangle\nonumber\end{equation}

\noindent The irreducible $L_0-$submodules are the one-, two-, and three-dimensional subspaces noted respectively above.   In particular, none of $X_{2, \,  j+2, \, 10}$, $X_{2, \,  j+2, \, 11}$,  or $X_{2, \,  j+2, \, 12}$ is $L_0-$isomorphic to either of the others.

	 In the case of $ X_{2, \,  j+2, \, 10},$ we have $(\ad [1, \, x])\cdot ((x^{(2)})^*\otimes \lambda + x^*\otimes \mu) = [1, \, \mu] = d,$ and $(\ad [x, \, x^{(2)}])\cdot d = [xy^{(j)}, \, x^{(2)}] = 2x^{(2)}y^{(j-1)} \neq 0,$ so by Lemma \ref{Lem:m2m2jp20}, 

\begin{equation}\label{eq:x23} X_{2, \,  j+2, \, 10} \cap L_{j+2} = 0.\end{equation}  

	In the case of $X_{2, \,  j+2, \, 11},$ we have $(\ad [1, \, x])\cdot (2\cdot x^*\otimes \kappa + 1^*\otimes \lambda) = 2\cdot [1, \, \kappa] + [\lambda, \, x] = 2e,$ and $2\cdot (\ad [1, \, x])\cdot e = (\ad [1, \, x])\cdot ((x^{(2)})^*\otimes xy^{(j)} + 2\cdot x^*\otimes y^{(j)}) = 2\cdot [1, \, y^{(j)}] = x^{(2)}y^{(j-1)} \neq 0,$ so by Lemma \ref{Lem:m2m2jp20}, 

\begin{equation}\label{eq:x24} X_{2, \,  j+2, \, 11} \cap L_{j+2} = 0.\end{equation} 

	In the case of $ X_{2, \,  j+2, \, 12},$ we have $(\ad [x, \, x^{(2)}])\cdot (x^*\otimes \lambda) = [\lambda, \, x^{(2)}] \equiv  d,$ and, as above, $(\ad [x, \, x^{(2)}])\cdot d = 2x^{(2)}y^{(j-1)} \neq 0,$  so by Lemma \ref{Lem:m2m2jp20}, 

\begin{equation}\label{eq:x25} X_{2, \,  j+2, \, 12} \cap L_{j+2} = 0.\end{equation}

\noindent  Thus, by transitivity (C),

\begin{equation}\label{eq:x2je14}X_{2, \,  j+1, \, 4}\cap L_{j+1} = 0.\end{equation} 

\noindent It follows that since $X_{2, \,  j+1, \, 4}$ is the irreducible $L_0-$submodule of the indecomposable $L_0-$module $X_{2, \,  j+1, \, 6},$ we must also have

\begin{equation}\label{eq:x2je16}X_{2, \,  j+1, \, 6}\cap L_{j+1} = 0.\end{equation}

Focusing next on $X_{2, \,  j+1, \, 5},$ we set

\begin{eqnarray}\nu &\eqdef&  1^* \otimes d \cr
\xi &\eqdef& 2\cdot x^*\otimes d + 1^* \otimes e \cr
\pi &\eqdef&  (x^{(2)})^*\otimes d + x^* \otimes e + 1^* \otimes f\nonumber\end{eqnarray}

\noindent Then $y \cdot \pi = 0,$ and $L_{-1}^*\otimes X_{2, \,  j+1, \, 5}$ contains the following  (non-isomorphic) $L_0-$submodules:

\begin{equation} X_{2, \,  j+2, \, 13} \eqdef \langle (x^{(2)})^* \otimes \nu + x^*\otimes \xi +  1^*\otimes \pi, \,  2\cdot x^*\otimes \nu + 1^*\otimes \xi, \, 1^*\otimes \nu\rangle\nonumber\end{equation}

\noindent and

\begin{equation}X_{2, \,  j+2, \, 14} \eqdef \langle (x^{(2)})^*\otimes \xi + x^*\otimes \pi + 1^*\otimes \nu, \, 2\cdot (x^{(2)})^*\otimes \nu + x^*\otimes \xi\rangle+\langle x^*\otimes \nu\rangle\nonumber\end{equation}

\noindent Furthermore,  $L_{-1}^*\otimes X_{2, \,  j+1, \, 5} = X_{2, \,  j+2, \, 13}  + X_{2, \,  j+2, \, 15},$ where modulo $X_{2, \,  j+2, \, 14}$, $X_{2, \,  j+2, \, 15}$  is spanned by $(x^{(2)})^*\otimes \nu$,  $2\cdot 1^*\otimes \nu + (x^{(2)})^*\otimes \xi,$ and $x^*\otimes \nu + 1^*\otimes \xi + (x^{(2)})^*\otimes \pi.$ We have

\begin{eqnarray} (\ad [x^{(2)}, \, 1])\cdot ((x^{(2)})^*\otimes \nu) &=& 2d\cr
(\ad [x, \, x^{(2)}])\cdot (2\cdot 1^*\otimes \nu + (x^{(2)})^*\otimes \xi) &=& 2d\cr
(\ad [(1, \, x])\cdot (x^*\otimes \nu + 1^*\otimes \xi + (x^{(2)})^*\otimes \pi) &=& 2d\nonumber\end{eqnarray}

\noindent so, since from above $(\ad [x, \, x^{(2)}])\cdot d = 2x^{(2)}y^{(j-1)} \neq 0,$ it follows from  Lemma \ref{Lem:m2m2jp20} that 

\begin{equation}\label{eq:s28}X_{2, \,  j+2, \, 15}\cap L_{j+2} \subseteq X_{2, \,  j+2, \, 14}\cap L_{j+2}\end{equation}

	Focusing on $X_{2, \,  j+2, \, 14},$ we have $(\ad [1, \, x])\cdot (2\cdot (x^{(2)})^*\otimes \nu+x^*\otimes \xi) = [1, \, \xi] = e,$ and $(\ad[1, \, x])\cdot e = [1,\, y^{(j)}] = 2x^{(2)}y^{(j-1)} \neq 0,$  so by Lemma \ref{Lem:m2m2jp20}, 

\begin{equation}\label{eq:x27} X_{2, \,  j+2, \, 14} \cap L_{j+2} = 0.\end{equation}

\noindent   It follows from \eqref{eq:s28} and \eqref{eq:x27} that  

\begin{equation}\label{eq:x28} X_{2, \,  j+2, \, 15}\cap L_{j+2} = 0.\end{equation}

	Note that  

\begin{equation}\label{eq:z26}[L_{-2}, \, X_{2, \,  j+2, \, 13}] = 0.\end{equation}

	We have observed before that $L_{j+2}$ is contained in $L_{-1}^*\otimes L_{-1}^*\otimes L_{-1}^*\otimes S_{j-1},$ whose intersection with $L_{j+2}$ has, in view of \eqref{eq:z210}, \eqref{eq:z218}, \eqref{eq:z216}, and \eqref{eq:z26}, zero bracket with $L_{-2}.$  In other words,
$[L_{-2}, \, L_{j+2}] = 0.$  But this contradicts Lemma \ref{Lem:m2jp2ne0}, and shows that $j \not \equiv 2$ (mod 3).

	When $j \equiv 0 \, (\hbox{mod 3)},$ there is only one irreducible three-dimensional $L_0$-submodule $Q_j$ of $(L_{-1})^*\otimes S_{j-1}$.  It is spanned by

\begin{eqnarray} &x^* \otimes x^{(2)}y^{(j)}\cr
&(x^{(2)})^* \otimes x^{(2)}y^{(j)} + 2\cdot x^* \otimes xy^{(j)} \hbox{ and }\cr
&1^* \otimes x^{(2)}y^{(j)} + (x^{(2)})^* \otimes xy^{(j)} + x^* \otimes y^{(j)}\nonumber\end{eqnarray}

\noindent and is the only $L_0-$submodule of $(L_{-1})^*\otimes S_{j-1}$ which has zero bracket with $L_{-2}.$   By what we observed above after \eqref{eq:up}, $Q_j \subseteq S_j,$  Then $Q_j$ must equal $ H(2: \mathbf{n}, \omega)_{j}$ and the analysis of $L_{-1}^*\otimes Q_j$ would be the same as that of $L_{-1}^*\otimes S_{j-1}$ when $j \equiv 1$ (mod 3). In particular, $S_{j+1} \cap L_{-1}^*\otimes Q_j \neq 0.$  If $S_{j+1}$ were to contain an $L_0-$submodule $\Gamma_{j+1}$ corresponding to $X_{1, \, j, \, 1}$ (which has non-zero bracket with $L_{-2}$), we could consider the Lie algebra generated by $L_{\leq 0} \oplus S_1 \oplus \cdots \oplus S_{j-1} \oplus Q_j \oplus \Gamma_{j+1}$ and obtain a Lie algebra in which the minimal gradation degree $i > 2$ such that $[L_{-2}, \, S_i] \neq 0$ is greater than $j,$ as required.

	We can therefore assume that $S_{j+1}$ contains an $L_0-$submodule $Q_{j+1}$ corresponding to $X_{1,j,2}$ (which has zero bracket with $L_{-2}$), so that $S_{j+2} \cap L_{-1}^*\otimes Q_{j+1} \neq 0.$ 
If $S_{j+2}$ were to contain an $L_0-$submodule $\Delta_{j+2}$ corresponding to  $X_{1, \, j+1, \, 1}$ or $X_{1, \, j+1, \, 3}$ (both of which have non-zero bracket with $L_{-2}$), we can consider the Lie algebra generated by $L_{\leq 0} \oplus S_1 \oplus \cdots \oplus S_{j-1} \oplus Q_j \oplus Q_{j+1} \oplus \Delta_{j+2}$ and obtain a Lie algebra in which the minimal gradation degree $i > 2$ such that $[L_{-2}, \, S_i] \neq 0$ is greater than $j,$ as required.  We can therefore assume that $S_{j+2}$ contains an $L_0-$submodule $Q_{j+2}$ corresponding to $X_{1, \, j+1, \, 2}$ (which has zero bracket with $L_{-2}$), so that $S_{j+3} \cap L_{-1}^*\otimes Q_{j+2} \neq 0,$ etc.

	In view of Lemma \ref{Lem:ad1cubed} (which shows that $Q_{j+2} \cong Q_{j-1}$), we can repeat this process to either obtain a Lie in which the minimal gradation degree $i > 2$ such that $[L_{-2}, \, S_i] \neq 0$ is successively greater, or to eventually find that the highest gradation space $S_s$ of $S$ contains an $L_0-$submodule $Q_s$ which has zero bracket with $L_{-2}.$  However, by Lemma \ref{Lem:2.3}, $S_s$ is irreducible as an $L_0-$module, so we would have that $[L_{-2}, \, S_s] = 0.$  Since $S$ is generated by $S_s$ and $L_{-1},$ it would follow from the construction of $Q_s (= S_s)$ (See also Lemma \ref{Lem:2.3}.)  that $[L_{-2}, \, S_{>0}] = 0,$ to contradict, for example, \eqref{eq:6.00001}.

	Whatever the residue of $j$ modulo three, then, we have shown that $L$ can be assumed to contain a Lie subalgebra such that  if $i \geqq 2$ is minimal such that $L_{-2}$ has non-zero bracket with the $i^{th}$ gradation space of that subalgebra, then $i$  is greater than $j.$

	We can therefore conclude the following:  Let $r' = r$ if there have been no replacements of $L$, or let $r'$ be the height of the most recent  replacement of $L$ otherwise.  Then, in either case,  if $j \geqq 2$ is  minimal such that $[L_{-2}, \, S_j] \neq 0,$ and $j < r' - 2,$  we can find a subalgebra of $L$ such that the smallest $i \geqq 2,$ such that $L_{-2}$ has non-zero bracket with the $i^{th}$ gradation space of that subalgbra, is greater than $j.$  Thus, by induction, $L$ contains a subalgebra such that the smallest $i \geqq 2$ such that $L_{-2}$ has non-zero bracket with the $i^{th}$ gradation space of the subalgbra, is at least $r' - 2.$  Thus, by Lemma \ref{Lem:6.02}, $(j = ) \,k \leqq 2,$ and we can now conclude that in any case

\begin{equation}j\leqq 2.\nonumber\end{equation}

	To show that $j = 1,$  we will, for a contradiction, assume that $j = 2.$  We begin by using an inductive argument from \cite[Lemma 2.14]{Sk2} to show that the centralizer of $S_s$ in $L_{<0}$ is zero.  Denote by $Z$ the centralizer of $S_s$ in $L$.  Then $Z = \oplus Z_i$ is a homogeneous subspace of $L.$  Since $S_s$ is stable under $\ad L_{\geqq 0}$, $Z$ is, as well.  The component $Z_{-1}$ is an $L_0-$submodule of $L_{-1},$ and $Z_{-1}$ by definition has zero bracket with $S_s.$   Thus, by (B) and (C),

\begin{equation}\label{eq:6.3.19}Z_{-1} = 0\end{equation}

	We will show that (See \cite[Lemma 2.14]{Sk2}.)

\begin{equation}\label{eq:6.3.2}Z_i = 0\end{equation}

\noindent  for all $i <0$, proceeding by (downward) induction on $i.$  
Assume that $i < -1,$ and that
 
\begin{equation}\label{eq:6.3.205}Z_{\ell} = 0, \, i < \ell < 0.\end{equation}
   
\noindent   By analogy with the previous argument, we here (where $j = 2$) assume that $L$ is generated by  $L_{<0} + L_0 + S_1 + S_2.$  If we show that $[S_1, \, Z_i] = [S_2, \, Z_i]  = 0,$  then $[{L}_{>0}, Z_i] = 0,$ so that

\begin{equation}\sum_{\ell \geqq 0}(\ad L_{-1})^{\ell} Z_i\nonumber\end{equation}

\noindent would be an ideal of ${L}$ properly contained in ${S}$, the simple ideal of ${L}$, to contradict the simplicity of ${S}.$  Consequently, to verify that $Z_i = 0$, we need only show that $[S_1, \, Z_i] = [S_2, \, Z_i] = 0.$

	We first show that $[S_1, \, Z_i] = 0.$ Indeed, when $i = -2,$ we have by \eqref{eq:6.1.3} that $[Z_{-2}, \, S_1] = 0 ,$  and when $i < -2,$ we have by \eqref{eq:6.3.205} that $[S_1, \, Z_i] \subseteq Z_{i+1} = 0.$  Similarly,  $[S_2, \, Z_i] \subseteq Z_{2+i} = 0$ when $i < -2.$   If $i = -2$ and $[S_2, \, Z_{-2}] \neq 0,$ then since by Lemma \ref{Lem:6.001} $L_{-2}$ is an irreducible $L_0-$module, we would have by Lemma \ref{Lem:2.315} that  $0 \neq [S_2, \, Z_{-2}] = [S_2, \, L_{-2}] =  [S_1, \, L_{-1}] = L_0.$ (See  Corollary \ref{Cor:1.5}.)  But then we would have $[L_0, \, S_s] = [[S_2, \, Z_{-2}], \, S_s] = [S_2, \,[ Z_{-2}, \, S_s]] = 0,$ to contradict  \cite[Lemma 2.13]{Sk2}. (See  also Lemma \ref{Lem:2.31}.)  This contradiction shows that $[Z_{-2}, \, S_2] = 0,$ and completes the verification of \eqref{eq:6.3.2}.

	Since we are assuming that $j  = 2$, we have

\begin{equation}\label{eq:6.3.21} [L_{-2}, \, S_2] \neq 0\end{equation}

\noindent and 

\begin{equation}\label{eq:6.3.22}[L_{-2}, \, S_1] = 0\end{equation}

\noindent We will first address the case in which $q > r,$ and then deal with the case in which $r \geqq q.$  Thus, assume first that

\begin{equation} q > r.\nonumber\end{equation}

\noindent Since (as we observed at the beginning of this section) $S_1$ is an irreducible $L_0-$module, it follows from \eqref{eq:6.3.2} that

\begin{equation}\label{eq:6.3.225} [L_{-r+1}, \, S_r] = S_1 \end{equation}

\noindent Then by \eqref{eq:6.3.22}, 

\begin{equation} [[L_{-r+1}, \, L_{-2}], \, S_r] = [L_{-r+1}, \, [L_{-2}, \, S_r]]\nonumber\end{equation}

\noindent If both sides of the equation were non-zero, then by irreducibility (B) they must both equal $L_{-1}.$  Then we would have, in view of Lemma \ref{Lem:2.31} and  \eqref{eq:6.3.2} and \eqref{eq:6.3.22}, that, if $r > 2,$

\begin{eqnarray} 0 \neq [[L_{-1}, \, S_1], \, S_r] &=& [[[L_{-r+1}, \, [L_{-2}, \, S_r]], \, S_1], \, S_r]\cr
 &=& [[[L_{-r+1}, \, S_1], \, [L_{-2}, \, S_r]], \, S_r]\cr
 &=&  [[L_{-r+1}, \, S_r], \, S_1], \, [L_{-2}, \, S_r]]\cr
 &\subseteq& [[S_1, \, S_1], \, [L_{-2}, \, S_r]]\nonumber \end{eqnarray}

\noindent which would imply that $[S_1, \, S_1] \neq 0,$ contrary to \eqref{eq:6.01}.  If, on the other hand, $r = 2,$ then the third line above becomes ``$\subseteq [[L_{-r+1}, \, S_r], \, S_1], \, [L_{-2}, \, S_r]] +  [[L_{-r+1}, \, S_r], \, S_1],$'' which leads to a similar contradiction. Thus, we can conclude that 

\begin{equation}\label{eq:6.3.23} [L_{-r+1}, \, [L_{-2}, \, S_r]] = 0\end{equation}

	 Now note that under the present assumptions, $r$ cannot equal four.  Indeed, by Lemma \ref{Lem:2.31},  $[[L_{-1}, \, S_1], \, S_r] =  [[L_{-1}, \, S_r], \, S_1]$  can be assumed to be non-zero and

\begin{equation}\label{eq:6.3.3}L_r \supseteq [S_1, \, S_{r-1}] \supseteq [S_1, \, [L_{-1}, \, S_r]]= [[L_{-1}, \, S_1], \, S_r] =  S_r\end{equation}

\noindent   If $r$ were four, we would have by \eqref{eq:6.3.3}, \eqref{eq:6.3.22}, and \eqref{eq:6.01} that

\begin{equation} [L_{-2}, \, S_r] =  [L_{-2}, \, S_4] =  [L_{-2}, \, [S_1, \, S_3]] = [S_1, \, S_1] = 0\nonumber\end{equation}

\noindent to contradict \eqref{eq:6.3.2}.  Similarly,  \eqref{eq:6.01} and  \eqref{eq:6.3.3} contradict one another when $r = 2$. Lastly, $r$ cannot equal three, either; indeed, by \eqref{eq:6.3.22}, $[[L_{-2}, \, L_{-2}], \, L_3] \subseteq [L_{-2}, \, L_1] = 0,$ so if $r$ were three, then by \eqref{eq:6.3.205}, $[L_{-2}, \, L_{-2}] = 0.$  Now suppose that $M_2$ is a non-zero  irreducible $L_0-$submodule of $L_2.$  Then by \eqref{eq:6.01}, $[L_{-1}, \, [L_1, \, M_2]] = [[L_{-1}, \, L_1], \, M_2],$ which is non-zero by \eqref{eq:nonullann}.  It follows that $[L_1, \, M_2]$ is a non-zero $L_0-$submodule of $S_3,$ which is $L_0-$irreducible by Lemma \ref{Lem:2.3}.  Then by \eqref{eq:6.3.22} and \eqref{eq:6.3.225}, we have 

\begin{equation} 0 \neq [L_{-2}, \, S_3] = [L_{-2}, \, [L_1, \, M_2]] = [L_1, \, [L_{-2}, \,  M_2]] \end{equation}

\noindent so that $[L_{-2}, \,  M_2] \neq 0.$  By \eqref{eq:nonullann} again, $[L_{-2}, \,  [L_{-2}, \,  M_2]] \neq 0.$  Consequently, in the Lie algebra $L_{-2} \oplus L_0 \oplus M_2,$ we have $B(L_{-2})_1 \neq 0.$  Then  (See Proposition \ref{Pro:1.4}.)  setting $L = B(L_{-2})$ and $W = L_1$ in the latter option of Lemma  \ref{Lem:1.6}, we arrive at a contradiction.  

	Since $r \geqq j = 2,$ it follows that

\begin{equation} r > 4\nonumber\end{equation}

	If we now bracket \eqref{eq:6.3.23} by $S_{r-2}$, we obtain

\begin{equation} 0 = [[L_{-r+1}, \, [L_{-2}, \, S_r]], \, S_{r-2}] = [[L_{-r+1}, \, S_{r-2}], \, [L_{-2}, \, S_r]]\nonumber\end{equation}

\noindent which (since in view of our assumption that $j=k=2$, we have $[L_{-2}, \, S_r] \neq 0$) would imply that 

\begin{equation}\label{eq:6.3.301}[L_{-r+1}, \, S_{r-2}] = 0\end{equation}

\noindent  since otherwise it would equal $L_{-1}$ by irreducibility (B), and transitivity (C) would be violated. 

	If we next bracket  \eqref{eq:6.3.23} by $S_{r-1},$ we get

\begin{equation}\label{eq:6.3.31}[[L_{-r+1}, \, S_{r-1}], \,[L_{-2}, \, S_r]] = 0\end{equation}  

	Now suppose that 

\begin{equation}\label{eq:6.315}[L_{-r}, \, S_{r-2}] = 0\end{equation}

\noindent  Then  we would have 

\begin{equation}0 = [[L_{-r}, \, S_{r-2}], \, S_{r-1}] = [[L_{-r}, \, S_{r-1}], S_{r-2}]\nonumber\end{equation}

\noindent which by (B) and (C) would imply that 

\begin{equation}\label{eq:6.3155} [L_{-r}, \, S_{r-1}] = 0\end{equation}

\noindent which in turn would (in view of (D)) imply that

\begin{equation}[L_{-r-1}, \, S_{r-1}] = [[L_{-r}, \, L_{-1}], \, S_{r-1}] \subseteq [L_{-r}, \, S_{r-2}]\nonumber\end{equation}

\noindent the right-hand side of which we have assumed to be zero.  But then we would have $[L_{-r-1}, \, S_{r-1}] = 0,$ so

\begin{equation} 0 = [0, \, S_r] =  [[L_{-r-1}, \, S_{r-1}], \, S_r] =  [[L_{-r-1}, \, S_{r}], \, S_{r-1}] = [L_{-1}, \,  S_{r-1}]\nonumber\end{equation}

\noindent by \eqref{eq:6.3.2} and irreducibility (B), to contradict transitivity (C).  We conclude that

\begin{equation}\label{6.3.32} [[L_{-r-1}, \, S_{r-1}], \, S_r] \neq 0\end{equation}

\noindent to imply that (See  \eqref{eq:6.3155}.)

\begin{equation} [L_{-r}, \, S_{r-1}] \neq 0
\nonumber\end{equation}

\noindent so that by irreducibility (B),

\begin{equation}\label{eq:6.3.4} [L_{-r}, \, S_{r-1}] = L_{-1}\end{equation}

\noindent  and to further imply that (See \eqref{eq:6.315}.)

\begin{equation}\label{eq:6.3.5} [L_{-r}, \, S_{r-2}] \neq 0\end{equation}

\noindent Then, since by Lemma \ref{Lem:6.001} $L_{-2}$ is an irreducible $L_0-$module, it follows that 

\begin{equation}\label{eq:6.3.5.5} [L_{-r}, \, S_{r-2}] = L_{-2}\end{equation}	

\noindent and, in view of \eqref{6.3.32}, that also

\begin{equation}\label{eq:6.3.8}[L_{-r-1}, \, S_{r-1}]  = L_{-2}\end{equation}

	If $[L_{-r}, \, [L_{-2}, \,  S_r]] = 0$, we would have by \eqref{eq:6.3.4} that 

\begin{eqnarray} 0 &=& [0, \, S_{r-1}]\cr 
&=& [[L_{-r}, \, [L_{-2}, \,  S_r]], \,  S_{r-1}]\cr
&=&  [[L_{-r}, \,  S_{r-1}], \, [L_{-2}, \,  S_r]\cr
&=&  [L_{-1}, \, [L_{-2}, \,  S_r]]\nonumber\end{eqnarray}

\noindent to contradict transitivity (C), since we are assuming that $(k = ) j = 2.$  Consequently,

\begin{equation}[L_{-r}, \, [L_{-2}, \,  S_r]] \neq 0\nonumber\end{equation}

\noindent so that, as above

\begin{equation}\label{eq:6.3.41}[L_{-r}, \, [L_{-2}, \,  S_r]] = L_{-2}\end{equation}

\noindent By \eqref{eq:6.3.8},

\begin{equation}\label{eq:6.3.9}[L_{-2}, \, S_2] = [[L_{-r-1}, \, S_{r-1}], \, S_2] = [[L_{-r-1}, \, S_2], \, S_{r-1}] \subseteq [L_{-r+1}, \, L_{r-1}]\end{equation}

	 Now, if $q > r+1,$ by \eqref{eq:6.3.2} and Lemma \ref{Lem:6.001}, we have $[L_{-r-2}, \, S_r] = L_{-2},$ so 

\begin{equation}\label{eq:6.3.92}[L_{-2}, \, S_2] \subseteq [[L_{-r-2}, \, S_r], \, S_2] =  [[L_{-r-2}, \, S_2], \, S_r] \subseteq [L_{-r}, \, S_r].\end{equation}

	If, on the other hand, $q = r+1,$ then, since by \eqref{eq:6.3.2} and irreducibility (B) $[L_{-r-1}, \, S_{r}] = L_{-1},$ we have

\begin{equation}[L_{-1}, \, S_1] = [[L_{-r-1}, \, S_{r}], \, L_1] \subseteq [L_{-r}, \, S_r] \subseteq [L_{-1}, \, S_1]\nonumber\end{equation}

\noindent since by (D), $[L_{-i},\, S_i] \subseteq [L_{-1}, \, S_1]$ for all $i$, $1 \leqq i \leqq \min\{q, \, r\},$ so that in particular, $[L_{-2}, \, S_2] \subseteq [L_{-1}, \, S_1] = [L_{-r}, \, S_r],$  so that \eqref{eq:6.3.92} holds when $q = r+1,$ also. Thus, we have by \eqref{eq:6.3.9} and \eqref{eq:6.3.92}, respectively,  (since we are assuming that $j = 2$)

\begin{eqnarray}\label{6.3.95}0 \neq  [L_{-2}, \, S_2] &\subseteq& [L_{-r+1}, \, S_{r-1}] \subseteq [L_{-1}, \, S_1]\cr
0 \neq [L_{-2}, \, S_2] &\subseteq& [L_{-r}, \, S_r] = [L_{-1}, \, S_1]
\end{eqnarray}

	By Lemma \ref{Lem:2.315},  we may assume that $[L_{-2}, \, S_2] = [L_{-1}, \, S_1]$ and thus conclude that

\begin{equation}\label{eq:6.3.955}   [L_{-r+1}, \, S_{r-1}] =  [L_{-r}, \, S_r]\end{equation}

 	Then we have from  \eqref{eq:6.3.23}, \eqref{eq:6.3.41} and \eqref{eq:6.3.955} that (See also \eqref{eq:6.3.31}.)

\begin{eqnarray} 0 &=& [0, S_{r-1}]\cr
&=&  [[L_{-r+1}, \, [ S_{r}, \, L_{-2}]], S_{r-1}]\cr
&=&  [[L_{-r+1}, \, S_{r-1}], [ S_{r}, \, L_{-2}]]\cr
&=&  [[L_{-r}, \, S_{r}], [ S_{r}, \, L_{-2}]]\cr
&=&  [[L_{-r}, \, [S_{r}, \, L_{-2}]], \, S_{r}]\cr
&=& [L_{-2}, \,  S_{r}]\cr
\nonumber\end{eqnarray}

\noindent to contradict \eqref{eq:6.3.2}.  This final contradiction enables us to conclude that if $q>r$,

\begin{equation}\label{eq:6.3.401} j = 1 \end{equation}

	Let us now assume that $q \leqq  r.$  For $i \geqq 1,$ set

\begin{equation}\label{eq:6.3.402} M_{-q+i} \eqdef (\ad S_1)^iL_{-q}\end{equation}

\noindent Then we would have, by \eqref{eq:6.3.22}, that $[L_{-2}, \,  M_{-q+i}] = 0$  for all $i.$  Furthermore, we would then have by induction that

\begin{equation}\label{eq:6.3.403} [L_{-i}, \,  M_{-q+i}] = 0, \hbox{ for all } i \geqq 2 \end{equation}

\noindent Now, if $[[L_{-i}, \, S_{q-1}], \, M_{-q+i}]$ were not equal to zero, then it would have to equal $L_{-1}$ by the irreducibility (B) of $L$, and we would have, by Lemmas \ref{Lem:2.1} and \ref{Lem:2.4} and \eqref{eq:6.3.403} that for $1 < i < q-1,$

\begin{eqnarray} L_{-q} &=& [L_{-q+1}, \, L_{-1}]\cr
&=&  [L_{-q+1}, \, [[L_{-i}, \, S_{q-1}], \, M_{-q+i}]]\cr
&=& [[L_{-i}, \, [L_{-q+1}, \, S_{q-1}]], \, M_{-q+i}]\cr
&\subseteq& [L_{-i}, \, M_{-q+i}]\ = 0\cr
\nonumber\end{eqnarray}

\noindent This contradiction shows that

\begin{equation}[[L_{-i}, \, S_{q-1}], \, M_{-q+i}] = 0, \, 2 \leqq i < q-1\nonumber\end{equation}

\noindent so that for $i < \frac{q}{2},$ we have by \eqref{eq:6.3.403} that

\begin{eqnarray} 0 &=& [L_{-q+i}, \, 0]\cr
&=& [L_{-q+i}, \, [[L_{-i}, \, S_{q-1}], \, M_{-q+i}]]\cr
&=& [[L_{-q+i}, \, [L_{-i}, \, S_{q-1}]], \, M_{-q+i}]\cr
\nonumber\end{eqnarray}

\noindent Since for $i = \frac{q}{2}$, we have $[L_{-q+i}, \, M_{-q+i}] = [L_{-\frac{q}{2}}, \, M_{-q+\frac{q}{2}}] = 0,$ also by \eqref{eq:6.3.403}, the above-displayed calculation is valid for $i \leqq  \frac{q}{2}.$   Consequently, if $[L_{-q+i},$ $[L_{-i}, \, S_{q-1}]]$ were not zero, then by the irreducibility (B) of $L$, it would equal $L_{-1}$, so that by Lemma 2.1, $M_{-q+i}$ would equal zero for $i \leqq \frac{q}{2}$.  In \cite{BGK} and \cite{GK1}, we proved the Main Theorem for all cases less than or equal to three.  Conseqently, we may assume that $q \geqq 4,$ so that  $M_{-q+i}$ would equal zero for all $i \leqq 2.$  However, when $i = 1,$  it follows from Lemma \ref{Lem:2.8} that $L_{-q+1} = [L_{-q}, \, S_1] = M_{-q+1}$,  so that $M_{-q+2}$ $=$ $[ [L_{-q}, \, S_1], \, S_1]$ $=$ $[L_{-q+1}, \, S_1].$    Moreover, it follows from Lemma \ref{Lem:1.6} (with $W = M_{-q+1}$ or $W = L_{-q}$) that $M_{-q+2} \neq 0$.  We conclude that $[L_{-q+i},$ $[L_{-i}, \, S_{q-1}]] = 0$ for $i \leqq 2.$  In particular,

\begin{equation}\label{eq:6.3.4055} [[L_{-q+2}, \, L_{-2}], \, S_{q-1}] = 0\end{equation}

\noindent Now, if $[L_{-q+2}, \, L_{-2}]$ were not equal to zero, it would, by Lemma \ref{Lem:2.4}, equal $L_{-q},$ and we would have $[L_{-q}, \, S_{q-1}] = 0,$ to contradict Lemma \ref{Lem:2.3}.  Thus, by Lemma \ref{Lem:2.8}, we must have

\begin{equation}0 = [L_{-q+2}, \, L_{-2}] \supseteq [[L_{-q}, \, S_2], \, L_{-2}] = [L_{-q}, \, [L_{-2}, \, S_2]]\nonumber\end{equation}

\noindent By Lemma \ref{Lem:2.315}, $ [L_{-2}, \, S_2] = [L_{-1}, \, S_1]$, so we have, in view of Lemmas \ref{Lem:2.8},

\begin{equation} 0 = [L_{-q}, \, [L_{-2}, \, S_2]]= [L_{-q}, \, [L_{-1}, \, S_1]] = [L_{-1}, \, [L_{-q}, \, S_1]] = [L_{-1}, \, L_{-q+1}]\nonumber\end{equation}

\noindent to contradict Lemma \ref{Lem:2.1}.  This contradiction shows that here, too, \eqref{eq:6.3.401} must be true.

	 Let $\tilde{L}$ be as in the statement of Lemma \ref{Lem:2.25}, and let $V_1$ be any irreducible $L_0-$submodule of $S_1.$ Because ${\tilde{L}}_{>0}$ is generated by $S_1$, we have $\{1\}-$transitivity (vi) in the negative part of
$\tilde{L}/M(\tilde{L}),$ and since by \eqref{eq:6.3.401} $[L_{-2}, \, S_1]$ $\neq 0,$ we can apply Lemma \ref{Lem:2.23} to $\tilde{L}/M(\tilde{L})$ to conclude that 

\begin{equation}\label{eq:m2v1nz} [L_{-2}, \, V_1] \neq 0\end{equation}

\noindent   Consequently, if $\tilde{\tilde{L}}$ is the Lie algebra generated by $L_{-1} \oplus L_0 \oplus V_1,$  then the depth $\tilde{\tilde{q}}$ of
$\tilde{\tilde{L}}/M(\tilde{\tilde{L}})$ (Again, see Theorem
\ref{Thm:1.3}.) is (also) greater than one.  Let $\tilde{\tilde{r}}$ denote the height of $\tilde{\tilde{L}}/M(\tilde{\tilde{L}}).$

 {\bf Case I: $\tilde{\tilde{q}} < \tilde{\tilde{r}}.$} Suppose first that the depth $\tilde{\tilde{q}}$ is less than $\tilde{\tilde{r}}.$ If $\tilde{\tilde{q}}$ is less than $q,$ then we can apply the
Main Theorem to $\tilde{\tilde{L}}/M(\tilde{\tilde{L}})$ to (inductively) conclude that the representation of $\tilde{\tilde{L}}_0 = L_0$ on $\tilde{\tilde{L}}_{-1} = L_{-1}$ is restricted, and see that the Main Theorem is true in this case. If $\tilde{\tilde{q}} = q,$ the Main Theorem follows from \cite{BGK}, \cite{GK1}, Sections 4 and 5, and Lemma \ref{Lem:3.1}.

 {\bf Case II: $\tilde{\tilde{q}} \geqq \tilde{\tilde{r}}.$} From now on, then, we will assume that $\tilde{\tilde{r}}$ is less than or equal to $\tilde{\tilde{q}}.$

 {\bf Case IIA: $\tilde{\tilde{q}} \geqq \tilde{\tilde{r}} > 1.$} Since $\tilde{\tilde{r}} > 1,$
it follows from the definition of $\tilde{\tilde{L}}$ that $[V_1, \, V_1] \neq 0,$
so that by Lemma \ref{Lem:2.28}, $\Ann_{L_0}V_1 = 0.$ Clearly, $B(V_1)$ (See Section 3.)
satisfies the conditions of the Main Theorem.  (Condition (E), for example, follows from the transitivity (C) of $L,$ which shows that actually $M(B(V_1)) = 0.$)  Now, $q$ is assumed to be greater than one, so we have
 
\begin{equation}0 \neq [L_{-1}, \, L_{-1}] = [(B(V_1))_1, \, (B(V_1))_1] \nonumber\end{equation}

\noindent so $B(V_1)$ is not degenerate.  Consequently, (since the depth $\tilde{\tilde{r}}$ of $B(V_1)$ is less than or equal to $\tilde{\tilde{q}}$ which is less than or equal to $q$) we can, as in Case I, apply the Main Theorem to conclude 
that the representation of $L_0'$ on $V_1$ is restricted, and that the representation of $L_0'$ on
$B(V_1)_1 = L_{-1}$ is restricted, as well (since $L_{-1} = B(V_1)_1 \subseteq \Hom(V_1, \, L_0);$ see also Lemma \ref{Lem:1.2}).

 {\bf Case IIB: $\tilde{\tilde{q}} \geqq \tilde{\tilde{r}} = 1.$} Since $\tilde{\tilde{r}} = 1,$ 
 
\begin{equation} [V_1, \, V_1] = 0.\nonumber\end{equation}

 By Corollary \ref{Cor:1.41}, either $L$ is degenerate and $r = 1,$ contrary to hypothesis, or $L$ and $\tilde{\tilde{L}}/M(\tilde{\tilde{L}})$ are  not  degenerate, in which case we can apply Proposition \ref{Pro:1.4} to $B((\tilde{\tilde{L}}/M(\tilde{\tilde{L}}))_1) = B(V_1)$ 
to conclude that $B(V_1)$ is isomorphic to $L(\epsilon)$ or 
$M,$ or is Hamiltonian (i.e., between $H(2:\mathbf{n},\omega)$ and $CH(2:\mathbf{n},\omega)$).
But in those cases, the one-component $(B(V_1))_{1} = L_{-1}$ is abelian; i.e., $[L_{-1}, \, L_{-1}] \subseteq M(\tilde{\tilde{L}}),$ so that by (D), $L_{-2} = [L_{-1}, \, L_{-1}] \subseteq M(\tilde{\tilde{L}}),$  so $[L_{-2}, \, V_1] = 0,$ to contradict \eqref{eq:m2v1nz}.  

 The proof of the Main Theorem is now complete.

\bigskip

Department of Mathematics, The Ohio State University at Mansfield,
Mansfield, Ohio  44906, USA

\bigskip

Department of Mathematics, Nizhni Novgorod State University,
Nizhni Novgorod 603600, Russia

\begin{thebibliography}{9}

\bibitem{BG} G.M. Benkart, T.B. Gregory, {Graded Lie
algebras with classical reductive null component}, Math. Ann. {\bf
285}, 1989, 85--98.
\bibitem{BGK} G.M. Benkart, T.B. Gregory, M.I.
Kuznetsov,  {On graded Lie algebras of characteristic three with
classical reductive null component},  The Monster and Lie Algebras,
Ohio State University Mathematical Research Institute Publications
{\bf 7}, 1998, 149--164.
\bibitem{BGP}  G.M. Benkart, T.B. Gregory, A.
Premet, {The Recognition Theorem for Graded Lie Algebras in Prime
Characteristic}, ${\mathcal Memoirs}$ of the American Mathematical Society, {\bf{197}}, no. 920 (second of 5 numbers) (2009).
\bibitem{BKK}     G.M. Benkart, A.I. Kostrikin, M.I.
Kuznetsov,  {The simple graded Lie algebras of characteristic
three with classical reductive component $L_0,$}  Comm. in Algebra
{\bf 24}, 1996, 223--234.
\bibitem{BW}     S. Berman, R.L. Wilson,
{Obstructions to modular classical simple Lie algebras},  Duke
Mathematical Journal {\bf 48}, 1981, 109--120.
\bibitem{B}      G. Brown,  {On the structure of some Lie
algebras of Kuznetsov},  Michigan Math. J. {\bf 39}, 1992, 85-90.
\bibitem{GK1}     T.B. Gregory, M.I. Kuznetsov, { On
depth-three graded Lie algebras of characteristic three with
classical reductive null component}, Communications in Algebra,
{\bf 33}, no. 9, pp. 3339-3371, 2004.
\bibitem{GK2}     T.B. Gregory, M.I. Kuznetsov,  
{Non-degenerate graded Lie algebras with a degenerate transitive subalgebra} (Russian), Contemporary Mathematics and its Applications, {\bf 60}, pp. 57-69, 2008.
\bibitem{J}    N. Jacobson,  {Lie Algebras},  Tracts
in Mathematics {\bf 10}, Interscience--New York, 1962.
\bibitem{K}    V.G. Kac, {The classification of simple
Lie algebras over a field of nonzero characteristic},  Izv. Akad.
Nauk SSSR, Ser. Mat. {\bf 34}, 1970, 385--408 (Russian), English
transl.  Math. USSR--Izv. {\bf 4}, 1970, 391--413.
\bibitem{KO}   A.I. Kostrikin, V.V. Ostrik,  {To the Recognition
Theorem for Lie algebras of characteristic 3},  Mat. Sbornik {\bf
186} (Russian), 1995, 73--88; translation in Sb. Math. {\bf
186}, 1995, no. 10, 1461--1475.
\bibitem{Sk1}   S.M. Skryabin, { New series of simple Lie
algebras of characteristic 3}, Mat. Sb. {\bf 183}, 1992, 3--22
(Russian), English transl. Russian Acad. Sci. Sb. Math {\bf 70},
1993, 389--406.
\bibitem{Sk2}   S.M. Skryabin, On the structure of the graded Lie algebra associated with a noncontractible filtration,  J. Algebra {\bf 197}, 1997, 178-230.
\bibitem{St}    H. Strade,  {Simple Lie
algebras over Fields of Positive Characteristic: I. Structure
Theory}, DeGruyter Expositions in Mathematics {\bf 38}, Walter de
Greyter--New York, 2004.
\bibitem{W}    B.J. Weisfeiler,  {On the structure of the
minimal ideal of some graded Lie algebras in characteristic $p >
0,$}  J. Algebra {\bf 53}, 1978, 344--361.

\end{thebibliography}
\end{document}